\documentclass[11pt]{article}

\usepackage[]{graphicx,float,latexsym,times}
\usepackage{amsfonts,amstext,amsmath,amssymb,amsthm}
\usepackage{caption2}

\usepackage{bm}
\usepackage{amssymb,mathtools} 
\usepackage{adjustbox}
\usepackage{multirow}

\long\def\symbolfootnote[#1]#2{\begingroup%
\def\thefootnote{\fnsymbol{footnote}}\footnote[#1]{#2}\endgroup}
\usepackage{titlesec} % Very fancy section title control: Sets sections as SQA requires

\titleformat{\section}{\large\bfseries}{\thesection.}{.5em}{}
\titlespacing*{\section}{0pt}{*3}{*2}
\titleformat{\subsection}{\normalfont\bfseries}{\thesubsection.}{.5em}{}
\titlespacing*{\subsection} {0pt}{*3}{*2}
\titleformat{\subsubsection}{\normalfont\bfseries}{\thesubsubsection.}{.5em}{}
\titlespacing*{\subsubsection} {0pt}{*3}{*2}

\theoremstyle{plain} %% italic text
\newtheorem{theorem}{Theorem}[section]
\newtheorem{lemma}{Lemma}[section]
\newtheorem{corollary}{Corollary}[section]

\theoremstyle{definition} %% or \theoremstyle{remark} will produce roman text

\newtheorem{remark}{Remark}[section]

\numberwithin{equation}{section} %% double numbering within sections

\usepackage[authoryear]{natbib}
\newcommand{\N}{ \mathbb{N} }

\newcommand{\wt}[1]{ \widetilde{ #1 } }

\newcommand{\eins}{{\bm 1}}
\newcommand{\matA}{{\bm A}}
\newcommand{\matB}{{\bm B}}

\newcommand{\matM}{{\bm M}}

\newcommand{\matS}{{\bm S}}
\newcommand{\matT}{{\bm T}}

\newcommand{\vecnull}{{\bm 0}}

\newcommand{\veck}{{\bm k}}

\newcommand{\vecu}{{\bm u}}
\newcommand{\vecs}{{\bm s}}

\newcommand{\vect}{{\bm t}}

\newcommand{\vecv}{{\bm v}}

\newcommand{\vecw}{{\bm w}}
\newcommand{\vecx}{{\bm x}}

\newcommand{\vecy}{{\bm y}}
\newcommand{\vecY}{{\bm Y}}

\newcommand{\bfsigma}{\bm \sigma}

\newcommand{\Var}{{\mbox{Var\,}}}
\newcommand{\Cov}{{\mbox{Cov\,}}}

\newcommand{\diag}{{\mbox{diag\,}}}

\newcommand{\cadlag}{c\`adl\`ag }

\begin{document}

\allowdisplaybreaks

\title{\textbf{\Large Detecting Changes in the Second Moment Structure of High-Dimensional Sensor-Type Data in a $K$-Sample Setting}}

\date{}

\maketitle

%%%%%%%%% Authors, affiliations %%%%%%%%%%%%%%%%%%%%%%%%%%

\author{
\begin{center}
\vskip -1cm

\textbf{\large Nils Mause and Ansgar Steland}

Institute of Statistics, RWTH Aachen University, \\
Aachen, Germany

\end{center}
}

\symbolfootnote[0]{\normalsize Address correspondence to Ansgar Steland,
	Institute of Statistics, RWTH Aachen University, W\"ullnerstra{\ss}e 3, 52062 Aachen, Germany; E-mail: steland@stochastik.rwth-aachen.de}

{\small \noindent\textbf{Abstract:} The $K$ sample problem for high-dimensional vector time series is studied, especially focusing on sensor data streams, in order to analyze the second moment structure and detect changes across samples and/or across variables cumulated sum (CUSUM) statistics of bilinear forms of the sample covariance matrix. In this model $K$ independent vector time series  $\vecY_{T,1},\dots,\vecY_{T,K}$  are observed over a time span $ [0,T] $, which may correspond to $K$ sensors (locations) yielding $d$-dimensional data as well as $K$ locations where $d$ sensors emit univariate data. Unequal sample sizes are considered as arising when the sampling rate of the sensors differs. We provide large sample approximations and two related change-point statistics, a sums of squares and a pooled variance statistic. The resulting procedures are investigated by simulations and illustrated by analyzing a real data set.}
\\ \\
%%%%%%%%% Key words %%%%%%%%%%%%%%%%%%%%%%%%%%
\\ \\

\section{INTRODUCTION}

In this paper we study the $K$ sample problem for high-dimensional data vectors assuming a time series framework and consider the problem to test for a change in the covariance matrix of the observed variables. This setting has many potential applications in big data analysis and data science in diverse areas such as internet data, finance, health science, climate research or in the field of production engineering used to illustrate the proposed methods by analyzing a real data set from the UCI Machine Learning Repository. Data sampled from sensor devices emitting data at certain sampling rates are pervasive in those areas and often lead to high-dimensional correlated time series thus motivating our work. Observing data from a large number of sensors located at $K$ locations easily fits into the studied framework, and the model assumptions as well as the sampling mechanism are somewhat tailored to this type of data. Nevertheless, the methods can be used to analyze quite general $K$-sample time series data.

Let us assume that a high-dimensional time series of such sensor-type data, $\{ \vecY_{t} \}$, is observed at,  say $ d = d_T$, sensors which are located at $ K \in \mathbb N $ distinct locations (sites). This formulation corresponds to a {\em single input single output} system, for example, when $ d_T $ sensors measure the same univariate quantity (such as speed, temperature or pressure), but our framework also applies to {\em single input multiple output} systems, where instead $K$ sensors at $K$ locations (defining the samples) generate a $ d_T $-dimensional correlated output signal.  The sensors measure the underlying quantity within a time interval $[0,T]$ for some $T>0$ at a site-specific sampling frequency $ \omega_j $, $ j = 1, \dots, K $, leading to sample sizes $ N_1, \dots, N_K $. Our asymptotics is for fixed sampling frequencies and $T \to \infty$ when $ N_j \sim N $ where $N$ is the total sample size. We consider a linear process time series framework for the data, which represents a natural model for sensor data, see, e.g., \cite{mari}, and, as shown in \cite{Steland2019}, covers a wide class of vector autoregressive moving average (VARMA) models with colored noise, amongst others. Linear processes are also a common framework in high-dimensional statistics, see, e.g., \cite{bstest}, \cite{WuHuangZheng2010} or \cite{WuMin2005}.

A common way to analyze high-dimensional data relies on projecting the given dataset onto a lower dimensional subspace and working with the projections instead of the high-dimensional vector itself. Well known examples for this procedure are principal component analysis (PCA), where one projects onto the eigenspaces associated to the covariance matrix $ \bm{\Sigma}_N $,  portfolio optimization, where the risk of the variance-minimal portfolio is given by such a projection statistic of the vector of asset returns, or data compression by projecting onto some known basis. 

When one is interested in drawing inference on the second moment structure of these projections, quadratic forms $\vecw' \bm{\Sigma}_N \vecw$ and more generally bilinear forms $\vecv' \bm{\Sigma}_N \vecw$, where $\vecv,\vecw \in \mathbb R^d$ are weighting vectors and $\bm{\Sigma}_N$ is the variance-covariance matrix of $\{\vecY_t\}$, appear naturally. Especially, if  $ \vecv = \vecw $ is an eigenvector, the  quadratic form is the associated eigenvalue leading to a method to draw inference on the eigenvalue structure. Other choices may correspond to variable selection, e.g., using projections derived from LASSO regressions.

The canonical nonparametric estimator for these bilinear forms is $\vecv' \hat{\bm{\Sigma}}_N \vecw$, where $\hat{\bm{\Sigma}}_N$ is the sample variance-covariance matrix. It is, however, well known that $ \hat{\bm{\Sigma}}_N $ is not a consistent estimator - and the same applies to the associated estimated eigenvalues and eigenvectors - when the dimension $d$ of the data is allowed to grow in an arbitrary way as the sample size $N$ increases, i.e. when usual assumptions such as $d \slash N \to y \in (0,1), \, N \to \infty,$ do not hold (cf. \cite{baiyin}, \cite{silverstein} or \cite{yin}), except  additional assumptions are made. For bilinear statistics of the form $\vecv' \bm{\Sigma}_N \vecw$ studied here, we also need to impose certain assumptions such as $ \ell_1 $-sparseness,  but there is no constraint on the growth of the dimension.

As estimated eigenvectors are an attractive choice for the projection vectors $ \vecv_N, \vecw_N $ when the problem of interest does not provide other guidance, let us briefly discuss some results on consistency of such estimators, although the literature on this subject is not yet matured and usually considers high-dimensional i.i.d. samples instead of time series. But these results indicate that $ \ell_q $-sparsity assumptions, $ q \le 2$, are quite common for high-dimensional data. For further discussion see \cite{Steland2019}. Classical PCA is known to be valid in spiked covariance models, cf. \cite{YataAoshima2013}, where the leading eigenvalues are larger than the remaining spectrum, and then can yield consistent estimates of the eigenvectors of $ \bm{\Sigma}_N $. For fixed eigenvalues \cite{BirnbaumJohnstonNadlerPaul2013} have shown consistency under $\ell_q$-constraints on the true eigenvectors. Further, under a joint $\ell_0 $-constraints on the rows of the matrix of eigenvectors consistent estimators have been studied by \cite{CaiMaWu2015}. In view of these results our assumptions on the projection vectors are not too restrictive, although the recent asymptotic results of \cite{WangFan2017} for estimated eigenstructures, which require diverging eigenvalues, may not be directly applicable. 

Focusing on the statistical problem to test for the presence of a change in the covariance structure, we study a CUSUM and a sum of squares (SSQ) statistic associated to a bilinear form of the sample covariance matrices. These statistics are related to bilinear forms of simple partial sums as arising when calculating the sample covariance matrix, such that large sample approximations of those partial sums yield the required approximations to construct statistical procedures. In \cite{steland1}, using martingale approximation techniques and building upon the results of \cite{kouritzin} and \cite{philipp}, appropriate  approximations have been derived. They hold without any constraint on the ratio of dimension and sample size as long as the weighting vectors $\vecv,\vecw$ are uniformly $\ell_1$ - bounded. For standardized versions using homogeneous estimators one may even use uniformly $ \ell_2$-bounded vectors. 

The present article contributes by extending these results to a $K$ - sample setting, studying CUSUM and SSQ change-point statistics in detail and elaborating on the pervasive case of sensor data. We derive  approximation results for the SSQ error processes and for bilinear forms based on the pooled sample variance-covariance matrix as an alternative approach. These approximations are then used to construct the two change-point tests, which aim at testing for the presence of a change (or structural break) in the covariance structure of sensor-type data in situations where a change may only affect a subset or all the given data, but there is no information which specific sample is affected. 

Finite sample properties in terms of size and the detection power, of both procedures are investigated by simulations. Further, the methods are illustrated by analyzing a real data set.

The paper is organized as follows. In Section 2 we will present the details of the used model and the necessary assumptions. Section 3  treats  weak approximation results for sum of squared error processes based on the original bilinear forms of \cite{steland2} and corresponding results for bilinear forms which are based on the pooled sample covariance matrix. The change-point tests are then introduced in Section 4. Results from extensive simulation studies are reported in Section 5. Section 6 illustrates the test by analyzing a real data set from production engineering. Proofs are provided in Section 7.

\section{MODEL AND MAIN ASSUMPTIONS}

To ease presentation, we stick to the single input single output sensor model. So let us assume that we have access to the information generated by $d_T \in \mathbb N$ sensors within each of $K \in \mathbb N$ different locations. At each of these $K$ locations the $d_T$ sensors quantitatively measure some underlying process $\{\epsilon_{j,k} :  k \in \mathbb Z\}, j=1,\dots,K,$ within the time interval $[0,T]$ for some $T > 0$, and then generate samples of output signals
\begin{align*} \vecY_{T,j,1},\dots,\vecY_{T,j,N_j}, \quad j=1,\dots,K,  \end{align*}
of sample size $N_j \in \mathbb N$ with a certain, possibly location-specific, sampling frequency $\omega_j \in (0,1]$. Then the output signals are $d_T$ - dimensional random vectors, i.e.
\begin{align*} \vecY_{T,j,i} = \left (Y_{T,j,i}^{(1)},\dots,Y_{T,j,i}^{(d_T)} \right )', \quad i=1,\dots,N_j, j=1,\dots,K,   \end{align*}
defined on probability spaces $(\Omega^{(j)},\mathcal{F}^{(j)},P^{(j)})$ for $j=1,\dots,K$, and as the sampling frequencies are allowed to be different for each location, the resulting sample sizes, given by $N_j = \lfloor \omega_j T \rfloor,j=1,\dots,K$, may also be different. %same probability space for all $T > 0$.
The output at time $i$ of the $\nu$th sensor (the $\nu$th coordinate of $ \vecY_{T,j,i}$) at location $j$ is assumed to follow a causal linear process (or filter) with real coefficients $ \{ c_{T,j,l}^{(\nu)} : l \geq 0 \} $, i.e.,
\begin{align} 
\label{ModelForY}
Y_{T,j,i}^{(\nu)} = \sum_{l=0}^\infty c_{T,j,l}^{(\nu)} \epsilon_{j,i-l}, \quad i=1,\dots,N_j, \nu=1,\dots,d_T, j = 1, \dots, K. \end{align}
As well known, causal linear filters are a standard model for sensor measurements often justified by the underlying circuits, see e.g., \cite{mari}. For results on time series models which can be represented resp. approximated by (\ref{ModelForY}), e.g. spiked covariance models and VARMA models, see \cite{Steland2019}.

We impose the following assumptions. 
The coefficients $\{ c_{T,j,l}^{(\nu)} : l \in \mathbb N_0 \}, \nu=1,\dots,d_T, j=1,\dots,K,$ in \eqref{ModelForY} are assumed to satisfy the decay condition
\begin{align} \label{decay}  \sup \limits_{T > 0} \max \limits_{1 \leq \nu \leq d_T} \vert c_{T,j,l}^{(\nu)} \vert^2 = \mathcal{O}(\min \{ l^{-\frac{3}{2} - \vartheta_j} , 1\}), \quad l \in \mathbb N_0, j=1,\dots,K,     \end{align}
for some $\vartheta_j \in \left (0 , \frac{1}{2} \right )$, and the underlying innovation processes $\{\epsilon_{j,k}\}_{k \in \mathbb Z}$ are supposed to be $K$ independent sequences of independent random variables with
\begin{align*} E(\epsilon_{j,k})=0, \quad E(\epsilon_{j,k}^2) = \sigma_{j,k}^2 \in (0,\infty) \qquad \text{and} \qquad \sup_{k \in \mathbb Z} E|\epsilon_{j,k}|^{4 + \delta_j} < \infty \end{align*}
for some $\delta_j > 0, j=1,\dots,K,$ such that we end up in a $K$ independent sample framework. With the assumption on the decay rate of the coefficients we not only cover short memory processes such as ARMA(p,q) - models, but also many long-range dependent series. We refer to \cite{steland1} and \cite{Steland2019} for discussions on this assumption. \\
\vspace{0.1cm}
\\Since we are interested in the second-moment structure of the high-dimensional output signals $\{\vecY_{T,j,i}\}$, we further assume $E(\vecY_{T,j,i}) = \vecnull$ for all $i \in \{1,\dots,N_j\}, j=1,\dots,K$, and define the $(d_T \times d_T)$ - dimensional sample variance-covariance matrices
\begin{align*}  \hat{\bm{\Sigma}}_{T}^{(j)} = \frac{1}{N_j} \sum_{i=1}^{N_j} \vecY_{T,j,i} \vecY_{T,j,i}', \quad j=1,\dots,K,  \end{align*}
which estimate the population variance-covariance matrices
\begin{align*}  \bm{\Sigma}_{T}^{(j)}  = E \left (\hat{\bm{\Sigma}}_{T}^{(j)} \right ) = \frac{1}{N_j} \sum_{i=1}^{N_j} E(\vecY_{T,j,i} \vecY_{T,j,i}') ,   \quad j=1,\dots,K,   \end{align*} 
and the pooled sample variance-covariance matrix 
\begin{align*}  \matS_T = \frac{1}{N} \sum_{j=1}^K \sum_{i=1}^{N_j} \vecY_{T,j,i} \vecY_{T,j,i}' =  \frac{1}{N} \sum_{j=1}^K N_j \hat{\bm{\Sigma}}_T^{(j)}, \end{align*}
where $N = \sum_{j'=1}^K N_{j'}$ is the total sample size. Note that these matrices depend on the time horizon $T$ through the sample sizes $N_j$ and that the following asymptotic results will hold for $T \to \infty$, which directly implies $N_j \to \infty$ for all $j=1,\dots,K.$

\section{ASYMPTOTICS}

The approximations are in terms of strong or weak approximations by Gaussian processes, which hold true for equivalent versions defined on a new probability space, on which the approximating process is also defined. Throughout the paper, we indicate those (equivalent) processes defined on the new space by a $ \tilde{\ } $. We shall, however, also describe how one may construct these approximations on the original probability space. 

\subsection{Approximations for the Sum of Squares Statistic}

As a preparation, let us introduce the following quantities.
For $k_j \in \mathbb N$ with $k_j \leq N_j$ define the $K$  $(d_T \times d_T)$ - dimensional matrix - valued partial sums 
\begin{align*} \hat{\bm{\Sigma}}_{T,k_j}^{(j)} = \sum_{i=1}^{k_j} \vecY_{T,j,i} \vecY_{T,j,i}', \qquad \text{and} \qquad \bm{\Sigma}_{T,k_j}^{(j)} = \sum_{i=1}^{k_j} E(\vecY_{T,j,i} \vecY_{T,j,i}') , \quad j=1,\dots,K.   \end{align*}
Further let $\hat{\bm \Sigma}_T^{(j)}(t_j) = N_j^{-1} \sum_{i=1}^{\lfloor t_j N_j \rfloor} {\vecY}_{T,j,i} {\vecY}_{T,j,i}'$ and $\bm \Sigma_{T}^{(j)}(t_j) = E(\hat{\bm \Sigma}_T^{(j)}(t_j)) $ for $t_j \in [0,1]$, such that $ \hat{\bm \Sigma}_T^{(j)}(1) $ is the sample covariance matrix of the $j$th sample. In a first step we shall study the bilinear form
\begin{align} \label{Abilinearform1} D_{T,k_j}^{(j)} = \vecv_{T,j}' \left ( \hat{\bm{\Sigma}}_{T,k_j}^{(j)} - \bm{\Sigma}_{T,k_j}^{(j)}  \right ) \vecw_{T,j}   , \quad k_j \geq 1, \end{align}
which appears naturally when analyzing the covariance of projections $\vecv_{T,j}' \vecY_{T,j,i}$ and $\vecw_{T,j}' \vecY_{T,j,i}$, where $\vecv_{T,j},\vecw_{T,j} \in \mathbb R^{d_T}$ are uniformly $\ell_1$ - bounded weighting vectors, i.e.
\begin{align*}   \sup \limits_{T > 0} \|\vecv_{T,j}\|_{\ell_1} < \infty, \qquad \sup \limits_{T > 0} \|\vecw_{T,j}\|_{\ell_1} < \infty.   \end{align*}
Note that \eqref{Abilinearform1} represents the error between the data driven bilinear form based on $\hat{\bm{\Sigma}}_T^{(j)}$ and the bilinear form based on the population covariance matrix $\bm{\Sigma}_T^{(j)}$. Further let
\begin{align} \label{Abilinearform2}  \mathcal{D}_{T}^{(j)}(t_j) &= N_j^{-1\slash 2} \vecv_{T,j}' \left ( \hat{\bm{\Sigma}}_{T,\lfloor t_j N_j\rfloor}^{(j)}  - \bm{\Sigma}_{T, \lfloor t_j N_j \rfloor}^{(j)} \right )  \vecw_{T,j} \\ & = \sqrt{N_j} \vecv_{T,j}' \left ( \hat{\bm{\Sigma}}_T^{(j)}(t_j) - \bm{\Sigma}_T^{(j)}(t_j)  \right ) \vecw_{T,j}, \quad t_j \in [0,1], T > 0, \nonumber  \end{align}
be the scaled \cadlag - process associated to \eqref{Abilinearform1}  and also define the bridge process
\begin{align} \label{Abridge} \Delta_{T}^{(j)}(t_j) = \mathcal{D}_{T}^{(j)} \left ( \frac{\lfloor t_j N_j \rfloor}{N_j} \right ) - \frac{\lfloor t_j N_j \rfloor}{N_j} \mathcal{D}_{T}^{(j)}(1), \quad t_j \in [0,1], T > 0.   \end{align}
For a stationary vector time series this process does not depend on the population variance-covariance matrix $\bm{\Sigma}_T^{(j)}$. \\

Let us introduce the long-run variance (LRV) parameters $\alpha_{j}^2$ 
\begin{align*} \alpha_{j}^2 =  \alpha_{j}^2(\{ \vecv_{T,j},\vecw_{T,j}\}) = \lim \limits_{T \to \infty} \Var \left ( \frac{1}{\sqrt{N_j}} \sum_{k_j=0}^{N_j} \eta_{T,j,k_j}   \right ), \quad j=1,\dots,K,  \end{align*}
associated to the univariate time series 
\begin{align*} \eta_{T,j,k_j} = \left ( \vecv_{T,j}' \vecY_{T,j,k_j}  \right ) \left ( \vecw_{T,j}' \vecY_{T,j,k_j}  \right ) - E \left ( \left ( \vecv_{T,j}' \vecY_{T,j,k_j}  \right ) \left ( \vecw_{T,j}' \vecY_{T,j,k_j}  \right )   \right ), k_j \geq 1, \end{align*}
for $j=1,\dots,K$. To exclude degenerate cases, we assume in the sequel that
\begin{align} \label{AssumptionPositiveAlphas} \alpha_j^2 > 0, \qquad j = 1, \dots, K.\end{align}

For the sum of squares associated to the bilinear forms, $  \sum_{j=1}^K  \left ( D_{T,k_j}^{(j)} \right )^2 $, we have the following approximation result. The approximation depends on the long-run-variance parameters $ \alpha_j^2 $ whose estimation will be discussed at the end of this section. 

\begin{theorem}
	\label{squaredtheorem}
	Let $(\vecv_{T,j},\vecw_{T,j}) \in \mathbb R^{d_T} \times \mathbb R^{d_T}$ be $K$ pairs of uniformly $\ell_1$ - bounded weighting vectors and  $\vecY_{T,j,i}$ be a vector time series following  model (\ref{ModelForY}). Then, for any $T >0$, there exist $K$ new probability spaces $(\tilde{\Omega}^{(j)}, \tilde{\mathcal{F}}^{(j)}, \tilde{P}^{(j)}), j=1,\dots,K$, on which one may define an equivalent version of $D_{T,k_j}^{(j)}$ and thus of $\mathcal{D}_{T}^{(j)}(t_j), t_j \geq 0$, denoted by $\tilde{D}_{T,k_j}^{(j)}$ and $\tilde{\mathcal{D}}_{T}^{(j)}(t_j)$, and independent standard Brownian motions $\{\tilde{B}_T^{(j)}(t_j) \, : \, t _j \geq 0 \}, j=1,\dots,K,$ which depend on the weighting vectors $(\vecv_{T,j},\vecw_{T,j})$, i.e. $\tilde{B}_T^{(j)}(t_j) = \tilde{B}_T^{(j)}(t_j,\vecv_{T,j},\vecw_{T,j}),$ 
	such that for some $\lambda_j > 0$ and constants $C_{T,j} < \infty, j=1,\dots,K,$ on the product space $(\tilde{\Omega},\tilde{\mathcal{F}},\tilde{P})$ of the $K$ probability spaces $(\tilde{\Omega}^{(j)}, \tilde{\mathcal{F}}^{(j)}, \tilde{P}^{(j)})$ it holds for all $k_j \in \mathbb N$ with $k_j \leq N_j$
	\begin{align} \label{squaredthmbound} \left \vert \sum_{j=1}^K  \left ( \tilde{D}_{T,k_j}^{(j)} \right )^2 - \sum_{j=1}^K \left ( \alpha_{T,j} \tilde{B}_T^{(j)} (k_j) \right )^2  \right \vert = \sum_{j=1}^K \alpha_{T,j} \mathcal{O}_P \left (\tilde{C}_{T,j} k_j^{1 - \lambda_j} \right ), \quad \tilde{P}-\text{a.s.},  \end{align}
	where $\tilde{C}_{T,j} = \max \{ C_{T,j}^2, C_{T,j}\}$. If additionally $\alpha_{T,j} \to \alpha_j^* \in (0,\infty)$ and $\tilde{C}_{T,j} N_j^{-\lambda_j} = o(1)$ for $j=1,\dots,K$ as $T \to \infty$, then we have the weak approximations 
	\begin{align} \label{squaredapprox1} & \sup \limits_{\vect \in [0,1]^K} \left \vert \sum_{j=1}^K  \left ( \tilde{\mathcal{D}}_{T}^{(j)}(t_j) \right )^2 - \sum_{j=1}^K \left ( \alpha_{T,j} \tilde{B}_T^{(j)} \left ( \frac{\lfloor t_j N_j \rfloor}{N_j}  \right ) \right )^2 \right \vert = o_{P}(1),  \quad T \to \infty,  \end{align} 
	and
	\begin{align} \label{squaredapprox2} &  \sup \limits_{\vect \in [0,1]^K} \left \vert \sum_{j=1}^K \left ( \tilde{\Delta}_T^{(j)}(t_j) \right )^2 - \sum_{j=1}^K \left [ \alpha_{T,j} \overline{\tilde{B}}_T^{(j)} \left ( \frac{\lfloor t_j N_j \rfloor}{N_j} \right )     \right ]^2  \right \vert  = o_{P}(1), \quad T \to \infty,  \end{align}
	where $\overline{\tilde{B}}_{T}^{(j)}(t_j) = \tilde{B}_{T}^{(j)}(t_j) - t_j \tilde{B}_{T}^{(j)}(1), t_j \in [0,1], j=1,\dots,K$, is a Brownian bridge.
\end{theorem}

The above approximations carry over to the standardized statistics. 

\begin{corollary}
	\label{squaredcorollary1}
	Given the assumptions of Theorem \ref{squaredtheorem}, such that in particular the weak approximations \eqref{squaredapprox1} and \eqref{squaredapprox2} hold. Then we also have 
	\begin{align} \label{squaredapproxst1} \sup \limits_{\vect \in [0,1]^K} \left \vert \sum_{j=1}^K  \left ( \frac{1}{\alpha_{T,j}} \tilde{\mathcal{D}}_{T}^{(j)}(t_j) \right )^2 - \sum_{j=1}^K \left ( \tilde{B}_T^{(j)} \left ( \frac{\lfloor t_j N_j \rfloor}{N_j}  \right ) \right )^2 \right \vert = o_P(1), \quad T \to \infty,  \end{align}
	and
	\begin{align} \label{squaredapproxst2}  \sup \limits_{\vect \in [0,1]^K} \left \vert \sum_{j=1}^K \left ( \frac{1}{\alpha_{T,j}} \tilde{\Delta}_T^{(j)}(t_j) \right )^2 - \sum_{j=1}^K \left [ \overline{\tilde{B}}_T^{(j)} \left ( \frac{\lfloor t_j N_j \rfloor}{N_j} \right )     \right ]^2  \right \vert  = o_P(1), \quad T \to \infty.  \end{align}
\end{corollary}
Here the approximating function of Brownian motions in \eqref{squaredapproxst1} satisfies
\begin{align*} \sum_{j=1}^K \left ( \tilde{B}_T^{(j)} \left ( \frac{\lfloor t_j N_j \rfloor}{N_j}  \right )   \right )^2 &  \stackrel{d}{=} \sum_{j=1}^K \xi_T^{(j)}, \qquad \xi_T^{(j)} \sim \Gamma \left (\frac12, 2 \frac{\lfloor t_j N_j \rfloor}{N_j} \right).  \end{align*} 
Explicit expressions for the pdf and cdf of $\sum_{j=1}^K \xi_T^{(j)}$ can be found in \cite{moschopoulos}.
However, in the special case $t_j=1$ and $N_j = \tilde{N} \in \mathbb N$ for all $j=1,\dots,K$ this simplifies to
\begin{align*} \sum_{j=1}^K \left ( \tilde{B}_T^{(j)} \left ( 1 \right )  \right )^2 \sim \chi^2(K).  \end{align*} 
For our change-point test in the next section we will use the following Corollary.
\begin{corollary}
	\label{squaredcusumcorollary}
	Suppose that the assumptions of Theorem \ref{squaredtheorem} hold, which especially imply the weak approximations \eqref{squaredapprox1} and \eqref{squaredapprox2}. Then we have
	\begin{align} \label{squaredcusum1} \left \vert \sum_{j=1}^K \max \limits_{k_j \leq N_j} \left ( \frac{1}{\alpha_{T,j}} \tilde{\mathcal{D}}_T^{(j)} \left ( \frac{k_j}{N_j}  \right )  \right )^2  - \sum_{j=1}^K \max \limits_{k_j \leq N_j} \left ( \tilde{B}_T^{(j)} \left ( \frac{k_j}{N_j} \right ) \right )^2 \right \vert = o_{P}(1), \end{align} 
	and
	\begin{align} \label{squaredcusum2} \left \vert \sum_{j=1}^K \max \limits_{k_j \leq N_j} \left ( \frac{1}{\alpha_{T,j}} \tilde{\Delta}_T^{(j)} \left ( \frac{k_j}{N_j}  \right )  \right )^2  - \sum_{j=1}^K \max \limits_{k_j \leq N_j} \left ( \overline{\tilde{B}}_T^{(j)} \left ( \frac{k_j}{N_j} \right ) \right )^2 \right \vert = o_{P}(1),  \end{align}	
	as $T \to \infty.$
\end{corollary} 
Since, for $j \in \{1,\dots,K\}$ and $y >0$, 
\begin{align*} & \tilde{P}^{(j)} \left ( \sup \limits_{t_j \in [0,1]} \left (\tilde{B}^{(j)} \left ( t_j \right ) \right )^2  \leq y \right )  = \tilde{P}^{(j)} \left ( \left ( \sup \limits_{t_j \in [0,1]} \left \vert \tilde{B}^{(j)}(t_j)  \right \vert \right )^2 \leq y \right ) \\ & = \tilde{P}^{(j)} \left (  \sup \limits_{t_j \in [0,1]}   \left \vert \tilde{B}^{(j)} \left ( t_j  \right ) \right \vert   \leq \sqrt{y} \right ) \end{align*}
(and similarly for $\overline{B}^{(j)}$), the distributions of the approximating function of Brownian motions in \eqref{squaredcusum1} and \eqref{squaredcusum2} are given through $K$ convolutions of the independent marginal distributions, which can be calculated using the reflection principle for Brownian motions in the case of \eqref{squaredcusum1} and the well known Kolmogorov distribution in the case of \eqref{squaredcusum2}. The following explicit formulas for the marginal distributions can be found in \cite{shorack}:
\begin{align*} \tilde{P}^{(j)} \left ( \sup \limits_{t_j \in [0,1]} \left \vert \tilde{B}^{(j)} \left ( t_j  \right ) \right \vert \leq \sqrt{y}  \right )  = \frac{4}{\pi} \sum_{l_j=0}^\infty \frac{(-1)^{l_j}}{2l_j+1} \exp \left ( - \frac{(2l_j+1)^2 \pi^2}{8y}   \right )  \end{align*}
\begin{align*} \tilde{P}^{(j)} \left ( \sup \limits_{t_j \in [0,1]} \left \vert \overline{\tilde{B}}^{(j)} \left ( t_j  \right ) \right \vert \leq \sqrt{y} \right ) = \sqrt{\frac{2\pi}{y}} \sum_{l_j=1}^\infty \exp \left (  - \frac{(2l_j-1)^2 \pi^2}{8y} \right )   \end{align*}
It remains to discuss whether the weak approximation results of Theorem \ref{squaredtheorem} also hold true on the original probability space, the product space $(\Omega,\mathcal{F},P)$ of the $K$ original probability spaces $(\Omega^{(j)},\mathcal{F}^{(j)},P^{(j)}),\, j=1,\dots,K$. This can be achieved as long as one can define a new uniformly distributed random variable.

\begin{theorem}
	\label{squaredapproxtheoremonoriginalprobabilityspace}
	Suppose that the assumptions of Theorem \ref{squaredtheorem} hold, which especially imply the weak approximations \eqref{squaredapprox1} and \eqref{squaredapprox2}. 
	If the underlying probability space $(\Omega, \mathcal{F},P)$ is rich enough to carry, in addition to the $K$ vector time series $\{ \vecY_{T,j,i} \, : \, 1 \leq i \leq N_j, N_j \geq 1   \}, j=1,\dots,K$, a uniform random variable $U$ that is independent of $ ( \mathcal{D}_T^{(j)}(\bullet) )_{j=1}^K $, 
	then the weak approximation results \eqref{squaredapprox1} and \eqref{squaredapprox2} in Theorem \ref{squaredtheorem} also hold on $(\Omega, \mathcal{F},P)$. This means, on $(\Omega, \mathcal{F},P)$ there exist $K$ Brownian motions $ \{ B_{T}^{(j)}(t_j) \, : \, t_j \geq 0  \}$ with
	$%\begin{align*} 
	\{ B_{T}^{(j)}(t_j) \, : \, t_j \geq 0  \} \stackrel{d}{=} \{ \tilde{B}_{T}^{(j)} (t_j) \, : \, t_j \geq 0 \} ,$ %\end{align*}
	such that
	\begin{align} \label{weakapproxonoriginal1} \sup \limits_{\vect \in [0,1]^K} \left \vert \sum_{j=1}^K  \left ( \mathcal{D}_T^{(j)}(t_j) \right )^2 - \sum_{j=1}^K \left ( \alpha_{T,j} B_{T}^{(j)} \left ( \frac{\lfloor t_j N_j \rfloor}{N_j}  \right ) \right )^2 \right \vert = o_P(1)   \end{align}
	and
	\begin{align} \label{weakapproxonoriginal2} \sup \limits_{\vect \in [0,1]^K} \left \vert \sum_{j=1}^K  \left ( \Delta_T^{(j)}(t_j) \right )^2 - \sum_{j=1}^K \left ( \alpha_{T,j} \overline{B}_{T}^{(j)} \left ( \frac{\lfloor t_j N_j \rfloor}{N_j}  \right ) \right )^2 \right \vert = o_P(1)   \end{align}
	as $T \to \infty$, where $\overline{B}_{T}^{(j)}(t_j) = B_{T}^{(j)}(t_j) - t_j B_{T}^{(j)}(1), t_j \in [0,1], j=1,\dots,K$, are  Brownian bridges.
\end{theorem}

Clearly, Theorem \ref{squaredapproxtheoremonoriginalprobabilityspace} also ensures that the approximations in Corollary \ref{squaredcorollary1} and \ref{squaredcusumcorollary}  hold for the original processes on $(\Omega,\mathcal{F},P)$. Especially,
\begin{align} \label{squaredcusumapproxonoriginal1} \left \vert \sum_{j=1}^K \max \limits_{k_j \leq N_j} \left ( \frac{1}{\alpha_{T,j}} \mathcal{D}_T^{(j)} \left ( \frac{k_j}{N_j}  \right )  \right )^2 - \sum_{j=1}^K \max \limits_{k_j \leq N_j} \left ( B_T^{(j)}\left ( \frac{k_j}{N_j} \right )  \right )^2        \right \vert  = o_P(1)   \end{align}
and
\begin{align} \label{squaredcusumapproxonoriginal2} \left \vert   \sum_{j=1}^K \max \limits_{k_j \leq N_j} \left ( \frac{1}{\alpha_{T,j}} \Delta_T^{(j)} \left ( \frac{k_j}{N_j}  \right )  \right )^2  -   \sum_{j=1}^K \max \limits_{k_j \leq N_j} \left ( \overline{B}_T^{(j)} \left ( \frac{k_j}{N_j} \right ) \right )^2 \right \vert = o_P(1).   \end{align}

\subsection{Approximations for the Pooled Sample Variance-Covariance Matrix}

Let us now study approximations for bilinear forms based on the pooled sample variance-covariance matrix $\matS_T$. For $k_j \in \mathbb N$ with $k_j \leq N_j$ let us define the $(d_T \times d_T)$ - dimensional matrix - valued partial sum
\begin{align*} \matS_{T,k_1,\dots,k_K} & = \sum_{j=1}^K \sum_{i=1}^{k_j} \vecY_{T,j,i} \vecY_{T,j,i}' = \sum_{j=1}^K \hat{\bm{\Sigma}}_{T,k_j}^{(j)}. \end{align*} 
Moreover we need to define the bilinear form
\begin{align} \label{bilinearform1} D_{T,k_1,\dots,k_K} = \vecv_T' (\matS_{T,k_1,\dots,k_K} - E(\matS_{T,k_1,\dots,k_K})) \vecw_T, \quad k_j \leq N_j,  \end{align}
and the associated \cadlag - process
\begin{align} \label{bilinearform2} \mathcal{D}_T(t_1,\dots,t_K) = \frac{1}{\sqrt{N}} \vecv_T' (\matS_{T,\lfloor t_1 N_1 \rfloor,\dots,\lfloor t_K N_K \rfloor}  - E(\matS_{T,\lfloor t_1 N_1 \rfloor, \dots, \lfloor t_K N_K \rfloor})  ) \vecw_T   \end{align}
for $t_1,\dots,t_K \in [0,1]$, where, again, $\vecv_T,\vecw_T \in \mathbb R^{d_T}$ are uniformly $\ell_1$ - bounded weighting vectors. Now, however, it is important to note that we take the same weighting vectors $\vecv_T$ and $\vecw_T$ for all of the $K$ locations as the pooled sample variance-covariance matrix $\matS_T$ can be represented as a convex combination of the $K$ sample variance-covariance matrices $\hat{\bm{\Sigma}}_T^{(j)}$. Therefore, by building bilinear forms of $\matS_T$ we extract the same information from all of the $K$ samples simultaneously.

\begin{theorem}
	\label{Pooledtheorem}
	Let $\vecv_T,\vecw_T \in \mathbb R^{d_T}$ be uniformly $\ell_1$ - bounded weighting vectors and  $\vecY_{T,j,i}$ be a vector time series following the given model structure of the previous section. Then, for any $T>0$, on a richer probability space, say $(\tilde{\Omega},\tilde{\mathcal{F}},\tilde{P})$, there exists an equivalent version of  $D_{T,k_1,\dots,k_K}$ and thus of $\mathcal{D}_{T}(t_1,\dots,t_K)$, $t_1\dots,t_K \geq 0$, denoted by $\tilde{D}_{T,k_1,\dots,k_K}$ and $\tilde{\mathcal{D}}_T(t_1,\dots,t_K)$, and $K$ independent Brownian motions $\{\tilde{B}_T^{(j)}(t_j) \, : \, t _j \geq 0 \}, j=1,\dots,K,$ which depend on the weighting vectors $(\vecv_T,\vecw_T)$, i.e. $\tilde{B}_T^{(j)}(t_j) = \tilde{B}_T^{(j)}(t_j,\vecv_T,\vecw_T),$ and have variances
	\begin{align*} E(\tilde{B}_T^{(j)}(t_j))^2 = \alpha_{T,j}^2 t_j, \qquad \alpha_{T,j}^2 = \alpha_{T,j}^2(\vecv_T,\vecw_T) > 0, \, t_j \geq 0,  \end{align*}
	such that for all $k_1,\dots,k_K \in \mathbb N$ with $k_j \leq N_j$ for all $j=1,\dots,K$ a.s.
	\begin{align} \label{pooledbound}  & \left \vert \tilde{D}_{T,k_1,\dots,k_K} -  \tilde{\mathcal{G}}_{T}(k_1,\dots,k_K)  \right \vert \leq \sum_{j=1}^K C_{T,j}k_j^{\frac{1}{2} - \lambda_j},   \end{align}
	where $C_{T,j} < \infty$ are constants, $\lambda_j > 0$ for $j=1,\dots,K$ and $\{\tilde{\mathcal{G}}_{T} (t_1,\dots,t_K) : t_1, \dots, t_K \geq 0  \}$ is a Gaussian process given by 
	\begin{align*} \tilde{\mathcal{G}}_{T}(t_1,\dots,t_K) = \sum_{j=1}^K \alpha_{T,j} \tilde{B}_{T}^{(j)}(t_j),  \end{align*}
	with mean zero and covariance function
	\begin{align*} \gamma_{\tilde{\mathcal{G}}_T}(\vecs,\vect) = \Cov(\tilde{\mathcal{G}}_T(s_1,\dots,s_K), \tilde{\mathcal{G}}_T(t_1,\dots,t_k)) = \sum_{j=1}^K \alpha_{T,j}^2 \min\{ s_j, t_j  \},  \end{align*}
	for $\vecs=(s_1,\dots,s_K)', \vect= (t_1,\dots,t_K)'$ with $s_j,t_j \geq 0$ for $j=1,\dots,K$. If additionally $ C_{T,j} N_j^{-\lambda_j} = o(1)$ and $\frac{N_j}{N} \to \kappa_j \in (0,1]$ for all $j=1,\dots,K$ as $ T \to \infty,$ then we also have
	\begin{align} \label{pooledapproxx} & \sup \limits_{\vect \in [0,1]^K} \left \vert \tilde{\mathcal{D}}_{T}(t_1,\dots,t_K) - \tilde{\mathcal{G}}_{T} \left ( \frac{\lfloor t_1 N_1 \rfloor}{N} , \dots, \frac{\lfloor t_K N_K  \rfloor}{N} \right )  \right \vert  = o(1),  \end{align}
	$\tilde{P}$-a.s., for $T \to \infty.$
\end{theorem}

\begin{remark}
	(i) The above strong approximation result \eqref{pooledapproxx} is based on $K$ independent approximations of the bilinear forms based on the sample variance-covariance matrix of each separate sample by their respective Brownian motion $\tilde{B}_T^{(j)}$, as shown by \cite{steland1}, and thus $(\tilde{\Omega},\tilde{\mathcal{F}},\tilde{P})$ is the product space of the probability spaces of the $K$ Brownian motions $\tilde{B}_T^{(j)},j=1,\dots,K.$ 
	
	(ii) The assumption of independent samples is required and cannot be weakened easily. It is, however, a routine assumption for $K$ sample problems.
\end{remark}

Next we present some special cases of the previous Theorem.

\begin{corollary}
	\label{pooledcorollary1}
	Let the conditions of Theorem \ref{Pooledtheorem} hold and suppose that additionally $\alpha_{T,j} \to \alpha_j^* \in (0,\infty)$ for all $j=1,\dots,K$ as $T \to \infty$. Then we have
	
	\begin{itemize}
		\item[(i)] the strong approximation
		\begin{align}  \label{StrongApproxSC} \sup \limits_{ \vect \in [0,1]^K}  \left \vert  \tilde{\mathcal{D}}_T^* (t_1,\dots,t_K) - \tilde{\mathcal{G}}_T^* \left ( \frac{\lfloor t_1 N_1 \rfloor}{N}, \dots, \frac{\lfloor t_K N_K \rfloor}{N} \right )   \right \vert = o(1),   \end{align}
		$\tilde{P}-$a.s., as $T \to \infty$ for the scaled process $\tilde{\mathcal{D}}_T^*(t_1,\dots,t_K) = (\sum_{j'=1}^K \alpha_{T,j'}^2 )^{-1 \slash 2} \tilde{\mathcal{D}}_T(t_1,\dots,t_K)$, \, $t_1,\dots,t_K \geq 0$, where
		\begin{align*} \tilde{\mathcal{G}}_T^*(u_1,\dots,u_K) = \frac{1}{\sqrt{\sum_{j'=1}^K \alpha_{T,j'}^2}}  \tilde{\mathcal{G}}_T \left ( u_1,\dots, u_K \right), \quad u_1,\dots,u_K \geq 0,   \end{align*}
		is again a Gaussian process with mean zero and covariance function
		\begin{align*} \gamma_{\tilde{\mathcal{G}}_T^*}(\vecu,\vecs) & = \Cov  ( \tilde{\mathcal{G}}_T^* (u_1,\dots,u_K)    , \tilde{\mathcal{G}}_T^*(s_1,\dots,s_K) ) = \frac{1}{\sum_{j'=1}^K  \alpha_{T,j'}^2} \sum_{j=1}^K \alpha_{T,j}^2 \min \{u_j, s_j\}      \end{align*}
		for $\vecu = (u_1,\dots,u_K)', \vecs = (s_1,\dots,s_K)'$ with $u_j,s_j \geq 0$ for $j=1,\dots,K$,
		
		\item[(ii)] the central limit theorem
		\begin{align} \label{Pooledzgws} \frac{N}{\sqrt{\sum_{j'=1}^K \alpha_{T,j'}^2 N_{j'}}} \vecv_T' \left ( \matS_T - E(\matS_T) \right ) \vecw_T  \stackrel{d}{\longrightarrow} \mathcal{N}(0,1), \qquad T \to \infty,   \end{align}
		for the original bilinear form defined on the product space $(\Omega, \mathcal{F},P)$ of the original $K$ probability spaces $(\Omega^{(j)},\mathcal{F}^{(j)},P^{(j)})$ based on the original processes $\{\vecY_{T,j,i} \, : \, i \geq 1\}$.

	\end{itemize}
\end{corollary}

The above results are usually not directly applicable, since they require us to know the covariance matrices $ \bm{\Sigma}_{T}^{(j)} = E (\hat{\tilde{\bm{\Sigma}}}_T^{(j)})  $, $ j = 1, \dots, K $. Having in mind the application to change-point testing, we replace them by the sample analogs $ N_j^{-1/2} \mathcal{D}_T^{(j)}(1)  $ leading to the the bridge process

\begin{align*} & {\Delta}_{T}(t_1,\dots,t_K) = {\mathcal{D}}_{T} \left ( t_1,\dots, t_K \right ) - \frac{1}{\sqrt{N}} \sum_{j=1}^K \frac{\lfloor t_j N_j \rfloor}{\sqrt{N_j}} {\mathcal{D}}_{T}^{(j)}(1),  \quad t_1,\dots,t_K \in [0,1].  \end{align*}

The following theorem yields an approximation by the corresponding Gaussian bridge process. It provides the approximation of the associated change-point test statistic 
under the null hypothesis of equal covariance matrices across samples, since then 
\[
\Delta_{T}(t_1,\dots,t_K) = \frac{1}{\sqrt{N}} \sum_{j=1}^K \vecv_T' \left (  \hat{\bm{\Sigma}}_{T,\lfloor t_j N_j \rfloor}^{(j)} - \lfloor t_j N_j \rfloor \hat{{\bm{\Sigma}}}_{T}^{(j)}  \right ) \vecw_T, \qquad T \ge 1.
\]

\begin{theorem}
	\label{Thmdtschlangeapprox}
	Given the assumptions and the constructions of Theorem \ref{Pooledtheorem}, especially the strong approximation \eqref{pooledapproxx}, we have
	\begin{align} \label{dtschlangeapprox} \sup \limits_{\vect \in [0,1]^K} \left \vert \tilde{\Delta}_{T}(t_1,\dots,t_K) -  \tilde{\Gamma}_{T} \left ( t_1,\dots,t_K   \right ) \right \vert  = o(1),   \end{align}
	$\tilde{P}$-a.s., for $T \to \infty,$ where
	\begin{align*} \tilde{\Gamma}_{T}(t_1,\dots,t_K) &= \tilde{\mathcal{G}}_{T} \left ( \frac{\lfloor t_1 N_1 \rfloor}{N}, \dots, \frac{\lfloor t_K N_K \rfloor}{N}   \right ) - \frac{1}{\sqrt{N}} \sum_{j=1}^K \frac{\lfloor t_j N_j \rfloor}{\sqrt{N_j}} \alpha_{T,j} \tilde{B}_{T}^{(j)}(1) \\ & \stackrel{d}{=} \sum_{j=1}^K \alpha_{T,j} \sqrt{\frac{N_j}{N}} \left [ \tilde{B}_T^{(j)} \left ( \frac{\lfloor t_j N_j \rfloor}{N_j} \right ) - \frac{\lfloor t_j N_j \rfloor}{N_j} \tilde{B}_T^{(j)}(1)     \right ]    . \end{align*}
\end{theorem}
Observe that here $\tilde{\Gamma}_T$ can be represented as a scaled sum of independent Brownian bridges.\\
\vspace{0.1cm}
\\Lastly, let us formulate the following CUSUM - Corollary which will be used in the next section in order to conduct our change-point test. We formulate it for the equivalent versions. 

\begin{corollary}
	\label{cusumcorollary2}
	Suppose that the assumptions of Theorem \ref{Pooledtheorem} hold implying the strong approximation \eqref{pooledapproxx}. Also let
	%\begin{align*}
	$ \mathcal{M}_T = \left \{(k_1,\dots,k_K) \in \mathbb N_0^K \, : \, k_1 \leq N_1, \dots, k_K \leq N_K   \right \}. $ 
	% \end{align*}
	Then we have
	\begin{align} \label{pooledcusum1}  & \left \vert   \max \limits_{\veck \in \mathcal{M}_T} \left \vert \tilde{\mathcal{D}}_{T} \left ( \frac{k_1}{N_1}, \dots, \frac{k_K}{N_K}   \right ) \right \vert - \max \limits_{\veck \in \mathcal{M}_T} \left \vert \tilde{\mathcal{G}}_{T} \left ( \frac{k_1}{N} , \dots, \frac{k_K}{N}   \right )   \right \vert \right \vert  = o(1)	 \end{align}
	and
	\begin{align} \label{pooledcusum2} \left \vert \max \limits_{\veck \in \mathcal{M}_T} \left \vert \tilde{\Delta}_{T} \left (\frac{k_1}{N_1},\dots,\frac{k_K}{N_K} \right ) \right \vert - \max \limits_{\veck \in \mathcal{M}_T} \left \vert \tilde{\Gamma}_{T} \left ( \frac{k_1}{N_1},\dots,\frac{k_K}{N_K} \right ) \right \vert \right \vert  = o(1),        \end{align}
	$\tilde{P}$-a.s., for $T \to \infty$.
\end{corollary}

Again, the above approximations can be constructed on the original probability space under weak additional assumptions. For brevity of presentation, we omit details and refer the reader to \cite{Mause2019}. 

\begin{theorem}
	\label{pooledapproxonoriginalspace}
	Suppose that the  assumptions of Theorem \ref{Pooledtheorem} hold. If the underlying probability space carries a uniform random variable $U$ that is independent from $ \matS_T $, then the approximations (\ref{pooledapproxx}), (\ref{StrongApproxSC}), (\ref{pooledcusum1}) and (\ref{pooledcusum2})  hold for the original processes in probability. 
\end{theorem}

\subsection{Estimation of the Unknown LRV Parameters}

In order to use the large approximations derived in the above results, we need estimators for the asymptotic variance parameters $ \alpha_j^2 $. Here we may apply the results obtained in \cite{steland1}: Let $m_{T,j}, T \geq 1$ be a sequence of lag truncation constants and $\{w_{m_{T,j},j,h} \, : \, h \in \mathbb Z, m_{T,j} \in \mathbb N, j=1,\dots,K \}$ be weights satisfying $w_{m_{T,j},j,h} \to 1$ as $m_{T,j} \to \infty$ for all $h \in \mathbb Z$ and $0 \leq w_{m_{T,j},j,h} \leq W_j$ for some constants $W_j, j=1,\dots,K$ and all  $m_{T,j} \geq 1$, $h \in \mathbb Z$. If additionally $\sup_{1 \leq \nu} \vert c_{T,j,l}^{(\nu)} \vert = \mathcal{O}(\max\{ l^{-1-\delta_j},1 \}), \{\epsilon_{j,k}\}_k$ is i.i.d. with $\max_k E \vert \epsilon_{j,k} \vert^8 < \infty$ and $m_{T,j} \to \infty$ with $m_{T,j}^2 \slash T \to 0$ as $T \to \infty$ then, in \cite{steland1} the $L_1$-consistency of the  following Bartlett type estimator 
\begin{align} \label{lrvestimate} \hat{\alpha}_{T,j}^2 =  \hat{\alpha}_{T,j}^2(\vecv_{T,j},\vecw_{T,j}) = \hat{\Gamma}_{T}^{(j)}(0) + 2 \sum_{h =1}^{m_{T,j}} w_{m_{T,j},h} \hat{\Gamma}_{T}^{(j)}(h),   \end{align}
for $\alpha_{j}^2$, has been shown, where for $\vert h \vert < N_j,$
\begin{align*} \hat{\Gamma}_{T}^{(j)}(h) = \frac{1}{N_j} \sum_{i=1}^{N_j - |h|} \left [ ( \vecv_{T,j}' \vecY_{T,j,i}) (\vecw_{T,j} ' \vecY_{T,j,i})  - \hat{\mu}_{T,j}  \right ] \left [ (\vecv_{T,j}' \vecY_{T,j,i+\vert h \vert})  (\vecw_{T,j}'  \vecY_{T,j,i+\vert h \vert} ) - \hat{\mu}_{T,j} \right ]   \end{align*}
and $\hat{\mu}_{T,j} = N_j^{-1} \sum_{r=1}^{N_j} (\vecv_{T,j}' \vecY_{T,j,r}) (\vecw_{T,j}' \vecY_{T,j,r})$ for $j=1,\dots,K$.\\
\vspace{0.1cm}

\section{CHANGE-POINT TESTS}

We aim at testing a-posteriori (i.e. off-line) whether there is a change in the covariances of at least one sample (location) assuming that all coordinates (sensors) are affected by the change at the same time. This situation occurs, for example, when considering a solar park with $K$ stations of sensors measuring the maximum power output of $d_T$ photovoltaic modules that they are attached to, and a change of weather or a failure of a solar inverter may affect only a part of the $K$ stations. This means, there is no information assumed about the specific location. 

The proposed test statistics are global tests in that they summarize the data from all samples. They also have power when only a fraction of the sensors resp. data vectors are affected by a change and/or are affected at different time points. Such cases can be treated by splitting samples and considering, in addition to the global tests, change-point tests and estimates as in \cite{steland1} and \cite{Steland2019} for each sample separately. For simplicity and brevity of presentation, we do not elaborate on these issues in greater detail. Further note that, within our framework, a change can be due to changes in the coefficients, changes in the variances of the innovations or both. However, the proposed test statistics do not require to estimate fit a time series model, which can be challenging and time-consuming, as all statistics are directly and easily computable from the data.

The  change-point testing problem of interest can be formulated as follows: We aim at testing the null hypothesis
\begin{align*} & H_0: \vecY_{T,j,1},\dots,\vecY_{T,j,N_j} \, \, \text{is a stationary time series with} \, \, E(\vecY_{T,j,i}) = \vecnull, \, \,  \\ & \qquad \Cov(\vecY_{T,j,i}) = \bm{\Sigma}_{T,0}^{(j)} \, \, \text{for all} \, \, i=1,\dots,N_j \, \, \text{and all locations} \, \,  j=1,\dots,K,  \end{align*}
against the alternative hypothesis
\begin{align*} H_1: \, & \text{There exists some $ \varnothing \neq \mathcal{J} \subseteq \{1,\dots,K\}$, such that for all}\, \, \tilde{j} \in \mathcal{J} \, \text{there exists some} \\ & \tau_{\tilde{j}}^* \in \{1,\dots,N_{\tilde{j}}\}, \, \, \text{such that} \, \, \vecY_{T,\tilde{j},1},\dots,\vecY_{T,{j},\tau_{\tilde{j}}^*} \, \, \text{is a}   \text{ stationary time series with} \\ & E(\vecY_{T,\tilde{j},i})=0 \, \, \text{and}   \, \, \Cov(\vecY_{T,\tilde{j},i}) = \bm{\Sigma}_{T,0}^{(\tilde{j})} \, \, \text{for all} \, \, 1 \leq i \leq \tau_{\tilde{j}}^* \\ &  \text{and} \, \,  \vecY_{T,\tilde{j},\tau_{\tilde{j}}^*+1},\dots,\vecY_{T,\tilde{j},N_j} \, \, \text{is a stationary time series with} \, \, E(\vecY_{T,\tilde{j},i})=0 \, \, \text{and} \\ &  \Cov(\vecY_{T,\tilde{j},i}) = \bm{\Sigma}_{T,1}^{(\tilde{j})} \, \, \text{for all} \, \, \tau_{\tilde{j}}^*+1 \leq i \leq N_{\tilde{j}}, \, \, \text{where} \, \, \Sigma_{T,0}^{(\tilde{j})} \neq \Sigma_{T,1}^{(\tilde{j})}.    \end{align*} 
Under the alternative, the set $\mathcal{J}$ is the subset of sample indices that are affected by the change, and for each $\tilde{j} \in \mathcal{J}$ this change occurs at index $\tau_{\tilde{j}}^*$ within the $\tilde{j}$-th data series $\vecY_{T,\tilde{j},1},\dots,\vecY_{T,\tilde{j},N_{\tilde{j}}}$, which corresponds to the (physical) time point $t_{\text{cp},\tilde{j}} =  \tau_{\tilde{j}}^* \slash \omega_{\tilde{j}}$ within the time interval $[0,T].$ 

Now let
\begin{align*} \vecv_T,\vecw_T \in \mathcal{W}_T & :=  \{ (\vecx_T,\vecy_T) \in \mathbb R^{d_T} \times \mathbb R^{d_T} \,  : \, \vecx_T' \bm{\Sigma}_{T,0}^{(j)} \vecy_T \neq \vecx_T' \bm{\Sigma}_{T,1}^{(j)} \vecy_T, \\ & \, \sup \limits_{T > 0} \|\vecx_T \|_{\ell_1} < \infty, \, \sup \limits_{T > 0} \|\vecy_T \|_{\ell_1} < \infty, \, \text{for some} \,\,   j \in \{ 1,\dots,K  \} \},  \end{align*} 
and note that under $H_1$ we have $\vert \mathcal{W}_T \vert > 0$.
For $\vecv_T,\vecw_T \in \mathcal{W}_T$  a change occurring in sample $\tilde{j} \in \mathcal{J}$ will also cause a change in the sequence of bilinear forms $\vecv_T' \bm{\Sigma}_{T}^{(\tilde{j})}(i) \vecw_T$, where $\bm{\Sigma}_T^{(\tilde{j})}(i) = \Cov(\vecY_{T,\tilde{j},i}), \, i \in \{1,\dots,N_{\tilde{j}} \}$, and thus $H_1$ implies
\begin{align*} H_1': \, & \text{There exists some $ \varnothing \neq \mathcal{J} \subseteq \{1,\dots,K\}$, such that for all} \, \, \tilde{j} \in \mathcal{J} \, \, \text{there exists some} \\ & \tau_{\tilde{j}}^* \in \{1,\dots,N_{\tilde{j}}\}, \, \, \text{such that} \, \, \sigma_{T,\tilde{j}}^2(i) = \sigma_{T,\tilde{j},0}^2 = \vecv_T' \bm{\Sigma}_{T,0}^{(\tilde{j})} \vecw_T  \, \, \text{for all} \, \,  1\leq i \leq \tau_{\tilde{j}}^* \, \, \text{and} \\ &  \sigma_{T,\tilde{j}}^2(i) = \sigma_{T,\tilde{j},1}^2 = \vecv_T' \bm{\Sigma}_{T,1}^{(\tilde{j})} \vecw_T \, \, \text{for all} \, \,  \tau_{\tilde{j}}^*+1 \leq i \leq N_{\tilde{j}}, \, \, \text{where} \, \, \sigma_{T,\tilde{j},0}^2 \neq \sigma_{T,\tilde{j},1}^2.  \end{align*} 

Assuming that the chosen projection vectors are contained in $ \mathcal{W}_T $ (if in doubt use, additionally, random projections), in view of the results of the previous section, we are now in a position to formulate change-point test statistics.  Firstly, when covariance matrices $\bm{\Sigma}_{T}^{(j)}$ are known one may use the sum of maximally selected  CUSUM of the SSQ statistics,
\begin{align*} Q_T = \sum_{j=1}^K \max \limits_{k_j \leq N_j} \left ( \frac{1}{\alpha_{T,j}} \mathcal{D}_T^{(j)} \left ( \frac{k_j}{N_j}  \right )  \right )^2 \stackrel{H_0}{=} \sum_{j=1}^K \max \limits_{k_j \leq N_j} \frac{1}{\alpha_{T,j}^2} \left [ \frac{1}{\sqrt{N_j}} \vecv_T' \left ( \hat{\bm{\Sigma}}_{T,k_j}^{(j)} - k_j \bm{\Sigma}_{T,0}^{(j)}  \right ) \vecw_T   \right ]^2.  \end{align*}
Alternatively, one may use the pooled sample variance-covariance matrix approach and consider the  maximized CUSUM test statistic
\begin{align*} V_{T} = \max \limits_{\veck \in \mathcal{M}_T}  \left \vert  \mathcal{D}_{T} \left ( \frac{k_1}{N_1}, \dots, \frac{k_K}{N_K} \right )    \right \vert  \stackrel{H_0}{=}  \max \limits_{\veck \in \mathcal{M}_T}  \left \vert \frac{1}{\sqrt{N}} \sum_{j=1}^K \vecv_T' \left ( \hat{\bm{\Sigma}}_{T,k_j}^{(j)} - k_j \bm{\Sigma}_{T,0}^{(j)}   \right ) \vecw_T   \right \vert .  \end{align*}
Secondly, if the variance-covariance matrices $\bm{\Sigma}_{T,0}^{(j)}$ are unknown, we can use  tests based on $\Delta_T$. Especially, we consider 
\begin{align*} \breve{Q}_T =  \sum_{j=1}^K \max \limits_{k_j \leq N_j} \left ( \frac{1}{\alpha_{T,j}} \Delta_T^{(j)} \left ( \frac{k_j}{N_j}  \right )  \right )^2  \stackrel{H_0}{=} \sum_{j=1}^K \max \limits_{k_j \leq N_j} \frac{1}{\alpha_{T,j}^2} \left [ \frac{1}{\sqrt{N_j}} \vecv_T' \left ( \hat{\bm{\Sigma}}_{T,k_j}^{(j)} - k_j \hat{\bm{\Sigma}}_{T}^{(j)}   \right ) \vecw_T   \right ]^2, \end{align*}
i.e. the sum of the maximally selected  $\Delta_T^{(j)}$s, and, as a test based on the pooled sample variance-covariance matrix, we use  
\begin{align*}  \breve{V}_{T} = \max \limits_{\veck \in M_T} \left \vert \Delta_{T} \left (\frac{k_1}{N_1},\dots,\frac{k_K}{N_K} \right ) \right \vert \stackrel{ H_0 }{=}  \max \limits_{\veck \in M_T} \left \vert \frac{1}{\sqrt{N}} \sum_{j=1}^K \vecv_T' \left (  \hat{\bm{\Sigma}}_{T,k_j}^{(j)} - k_j \hat{\bm{\Sigma}}_{T}^{(j)}  \right ) \vecw_T \right \vert.  \end{align*}

The following lemma summarizes the asymptotic null distributions of these tests.
 
\begin{lemma}
	\label{teststatisticslemma}
	Let the approximation results of Theorems \ref{squaredapproxtheoremonoriginalprobabilityspace} and $\ref{pooledapproxonoriginalspace}$ hold, such that the results of Corollary \ref{squaredcusumcorollary} as well as Corollary $3.6$ hold on the product space $(\Omega,\mathcal{F},P)$ of the $K$ original probability spaces $(\Omega^{(j)},\mathcal{F}^{(j)},P^{(j)}),j=1,\dots,K$. Then we have
	\begin{itemize}
		\item[(i)] $ Q_T \sim_{T \to \infty}  \sum_{j=1}^K \sup \limits_{s_j \in [0,1]} \left ( B^{(j)} \left ( s_j \right ) \right )^2$ \quad $ (ii) \ $ 	
		%	\item[(ii)] 
		$V_T \sim_{T \to \infty} \sup \limits_{\vecs \in  [0,1]^K } \left \vert  \sum_{j=1}^K \alpha_{j}^* \sqrt{\kappa_j} B^{(j)} \left ( s_j \right )   \right \vert$
		
		\item[(iii)] $\breve{Q}_T \sim_{T \to \infty} \sum_{j=1}^K \sup \limits_{s_j \in [0,1]} \left ( \overline{B}^{(j)} \left ( s_j \right ) \right )^2$ \quad $ (iv) $
		%	\item[(iv)] 
		$\breve{V}_T \sim_{T \to \infty}  \sup \limits_{\vecs \in [0,1]^K} \left \vert \sum_{j=1}^K \alpha_j^* \sqrt{\kappa_j} \overline{B}^{(j)}(s_j)       \right \vert$
	\end{itemize}	
\end{lemma} 

Critical values for the tests based on $V_T$ and $\breve{V}_T$ need to be simulated using estimates for the $\alpha_{T,j}$'s.

\section{SIMULATIONS}

We conducted a simulation study in order to examine the statistical properties of the proposed tests in terms of  level and power for a nominal significance level $ \alpha = 5\% $, focusing on the realistic case of an unknown covariance matrices $\bm{\Sigma}_{T,0}^{(j)}, j=1,\dots,K$. For $ K = 4 $ locations with $d \in \{10, 50, 200, 500, 1\,000\}$ sensors at each location autoregressive series were simulated over a time span of $  T = 1200 $ time points. With respect to the sample sizes four cases listed in Table \ref{samplesizetable} were investigated. 
\begin{table}[h]
	\begin{center}
		\begin{tabular}{|c|c|c|c|c|} \hline
			Case & $N_1$ & $N_2$ & $N_3$ & $N_4$ \\ \hline \hline
			I & 100 & 120 & 70 & 90 \\ 
			II & 300 & 250 & 350 & 180 \\
			III & 500 & 450 & 550 & 600 \\ 
			IV & 1000 & 900 & 1100 & 950 \\ \hline
		\end{tabular}
		\caption{Sample sizes used in the simulations}
		\label{samplesizetable}
	\end{center}
\end{table}

For those sample sizes AR(1) series were simulated with pre-change AR coefficients $\rho_{\nu,0} = 0.1 + 0.5 \nu \slash d$ and, when studying a change in the coefficients, after-change parameters  $\rho_{\nu,1} = 0.4 + 0.5 \nu \slash d$  for all $\nu =1,\dots,d, l=0,\dots,50$ and $j=1,\dots,4$. Having practical applications in mind, we stick to the case that all sensors at all locations follow the same model. The innovation processes were simulated as i.i.d. mean zero normal distributions with standard deviations given by $\tilde{\bfsigma}_{0} = (1,1.5,0.7,1)'$ if there is no change-point, and with after-change values $\tilde{\bfsigma}_{1} = (1,0.7,1.2,1)'$ when studying a change in the innovations.

The projection vectors $\vecw$ were drawn from a Dirichlet distribution of order $d$ with random parameters $\vartheta_1,\dots,\vartheta_d \sim \mathcal{U}[0,1],$ such that the $\ell_1$-condition for the weighting vectors holds. 

To estimate the long-run variance parameter one may either use a learning sample or use in-sample estimation. Both cases were investigated. The  lag-truncation constant and bandwidth were adaptively chosen as proposed by \cite{andrews} with weights given by a quadratic spectral kernel function.\\
\vspace{0.1cm}

Let us first discuss the results when using a learning sample  of size $L_1 = 500$ and $L_2 = 1\,000$, respectively. The simulated type I error rates for the sum of squared errors test  $\breve{Q}_T$ (for unknown covariances matrices under $H_0$) are provided in Table \ref{squaredtestunknowncovalpha}. Those for the pooled test are given in Table \ref{pooledtestunknowncovalpha}. Here, as in all reported simulations, each table entry is based on $10\,000$ runs. One can notice that the level is quite stable with respect to the dimension $d$, but both tests tend to overreact. This effect is, however, very mild for the pooled test, which is much more accurate in terms of the type I error rate. 

\begin{table}[h]
	\begin{center}
		\begin{adjustbox}{max width=\textwidth}
			\begin{tabular}{|c|c|c|c|c|c||c|c|c|c|c|c|c|c|c|c|} \hline
				& \multicolumn{5}{|c|}{$L_1 = 500$} & \multicolumn{5}{|c|}{$L_2 = 1\,000$} \\ \hline \hline
				Case	& d=10 & d=50 & d=200 & d=500  & d=1000 & d=10 & d=50 & d=200 & d=500  & d=1000 \\ \hline \hline
				I &   9.51 & 9.81 &  10.34 &  9.29  &  9.52 & 7.86 & 8.84 &   8.40 &   8.68  &   8.54 \\ \hline
				II &  7.97 &  8.80  &  8.68 &  8.85 &   9.17 & 8.16 &  7.68 &  7.26 &  7.44  &  6.83 \\ \hline
				III &  9.65 &  9.40 &  8.82 &  8.67 &   8.77 & 8.65  & 7.58 &  7.26 &  7.13  &  7.21 \\ \hline
				IV &  8.80 &  8.58 &  8.79 &  9.45  &  8.71 & 6.22  & 6.99 &   7.58 &  7.74  &  7.31 \\ \hline
			\end{tabular}
		\end{adjustbox}
		\caption{Estimated error of first order of the sum of squared errors test for learning sample sizes $L_1 = 500$ and $L_2=1\,000$ (values in \%)}
		\label{squaredtestunknowncovalpha}
	\end{center}
\end{table}

\begin{table}[h]
	\begin{center}
		\begin{adjustbox}{max width=\textwidth}
			\begin{tabular}{|c|c|c|c|c|c||c|c|c|c|c|c|c|c|c|c|} \hline
				& \multicolumn{5}{|c|}{$L_1 = 500$} & \multicolumn{5}{|c|}{$L_2 = 1\,000$} \\ \hline \hline
				Case	& d=10 & d=50 & d=200 & d=500  & d=1000 & d=10 & d=50 & d=200 & d=500  & d=1000 \\ \hline \hline 
				I &   6.21 & 6.51 &  6.36 &  6.72 &   6.76 &  5.60 &  6.00 &  5.59 &  5.82  &  5.87 \\ \hline
				II &  6.46 &  6.74 &   6.60 &  6.37 &   6.53 & 6.29 &  5.37 &  5.40 &  5.35 &   5.52 \\ \hline
				III & 6.93 &  6.48 &  6.70 &  6.09  &  6.69 &  5.52 & 5.87 &  5.83 &  5.67 &   5.22 \\ \hline
				IV  & 6.87 & 6.32 &  7.16 &   6.50  &  6.95 & 5.97 & 5.96 &  5.77  & 5.47 &   5.85 \\ \hline
			\end{tabular}
		\end{adjustbox}
		\caption{Estimated error of first order of the pooled test for learning sample sizes $L_1 = 500$ and $L_2 = 1\,000$ (values in \%)}
		\label{pooledtestunknowncovalpha}
	\end{center}
\end{table}

To investigate the power we simulated a change occurring either after $240$, $600$ or $960$ time instants, affecting all four locations either very early, in the middle or quite late within the time interval $[0,T]$. The first set of simulations examines a change in the standard deviations of the innovations from $\tilde{\bfsigma}_{0} = (1,1.5,0.7,1)'$ to $\tilde{\bfsigma}_{1} = (1,0.7,1.2,1)'$, whereas the second set investigates a change in the coefficients as described above.

The results for the first set (change in the innovations) are given in Tables \ref{squaredtestunknowncovpowercpinsds}  and  \ref{pooledtestunknowncovpowercpinsd}.
For both tests the detection power is better for a change in the middle of the sample than for an early or late change. Further, the SSQ  test has higher power than the pooled test.

\begin{table}[h]
	\begin{center}
		\begin{adjustbox}{max width=\textwidth}
			\begin{tabular}{|c|c|c|c|c|c|c||c|c|c|c|c|c|c|c|c|} \hline
				&		& \multicolumn{5}{|c|}{$L_1 = 500$} & \multicolumn{5}{|c|}{$L_2 = 1\,000$} \\ \hline \hline
				Change	&	Case	& d=10 & d=50 & d=200 & d=500  & d=1000 & d=10 & d=50 & d=200 & d=500  & d=1000 \\ \hline 
				\multirow{4}{*}{\shortstack{after 240 \\ time instants}} & I &  0.8234 & 0.8000 & 0.8166 & 0.8315 &  0.8350 & 0.8428 & 0.8491 & 0.8234 & 0.8258 &  0.8102 \\
				& II & 0.9997 & 0.9998 & 0.9997 & 0.9997 &  0.9999 & 0.9994 & 1.0000 & 0.9998 & 0.9999 &  0.9996   \\ 
				& III & 1.0000 & 1.0000 & 1.0000 & 1.0000 & 1.0000 & 1.0000 & 1.0000 & 1.0000 & 1.0000 &  1.0000 \\ 
				&	IV & 1.0000 & 1.0000 & 1.0000 & 1.0000 & 1.0000 &  1.0000 & 1.0000 & 1.0000 & 1.0000 &  1.0000 \\ \hline \hline
				\multirow{4}{*}{\shortstack{after 600 \\ time instants}} & I &  0.9602 & 0.9488 & 0.9398 & 0.9438 &  0.9415 & 0.9091 & 0.931 & 0.9397 & 0.9433 &  0.9437  \\
				& II &  1.0000 & 1.0000 & 1.0000 & 1.0000 &  1.0000 & 1.0000 & 1.0000 & 1.0000 & 1.0000 &  1.0000   \\ 
				& III & 1.0000 & 1.0000 & 1.0000 & 1.0000 &  1.0000 &  1.0000 & 1.0000 & 1.0000 & 1.0000 & 1.0000  \\ 
				& IV  & 1.0000 & 1.0000 & 1.0000 & 1.0000 &  1.0000 &  1.0000 & 1.0000 & 1.0000 & 1.0000 & 1.0000 \\ \hline \hline
				\multirow{4}{*}{\shortstack{after 960\\ time intants }} & I & 0.6304 & 0.6436 & 0.6224 & 0.6392 &  0.6290 &  0.6439 & 0.6323 & 0.6305 & 0.6311 &  0.6264 \\ 
				&	II & 0.9992 & 0.9946 & 0.9953 & 0.9967 &  0.9962 & 0.9961 & 0.9953 & 0.9971 & 0.9960 & 0.9969   \\ 
				&	III& 0.9999 & 1.0000 & 1.0000 & 1.0000  & 1.0000 & 0.9999 & 1.0000 & 1.0000 & 1.0000 & 1.0000 \\ 
				&	IV & 1.0000 & 1.0000 & 1.0000 & 1.0000 &  1.0000 & 1.0000 & 1.0000 & 1.0000 & 1.0000 & 1.0000 \\ \hline

			\end{tabular}
		\end{adjustbox}
		\caption{Simulated power of the sum of squared errors test for a change in the standard deviations of the innovation processes with learning sample sizes $L_1=500$ and $L_2=1\,000$.}
		\label{squaredtestunknowncovpowercpinsds}
	\end{center}
\end{table}

\begin{table}[h]
	\begin{center}
		\begin{adjustbox}{max width=\textwidth}
			\begin{tabular}{|c|c|c|c|c|c|c||c|c|c|c|c|c|c|c|c|} \hline
				&		& \multicolumn{5}{|c|}{$L_1 = 500$} & \multicolumn{5}{|c|}{$L_2 = 1\,000$} \\ \hline \hline
				Change	&	Case	& d=10 & d=50 & d=200 & d=500  & d=1000 & d=10 & d=50 & d=200 & d=500  & d=1000 \\ \hline 
				\multirow{4}{*}{\shortstack{after 240 \\ time instants}} & I & 0.1217 & 0.1401 & 0.1294 & 0.1410 &  0.1390 & 0.1104 & 0.1249 & 0.1285 & 0.1288 &  0.1240  \\
				& II & 0.3417 & 0.3635 & 0.3691 & 0.3791 &  0.3473 & 0.2384 & 0.3567 & 0.3545 & 0.3525 &  0.3480  \\ 
				& III & 0.6656 & 0.6257 & 0.6508 & 0.6661 &  0.6725 & 0.8057 & 0.7157 & 0.6842 & 0.6703 & 0.6476  \\ 
				&	IV & 0.9776 & 0.9868 & 0.9850 & 0.9841 &  0.9829 & 0.9855 & 0.9904 & 0.9858 & 0.9924 &  0.9909  \\ \hline \hline
				\multirow{4}{*}{\shortstack{after 600 \\ time instants}} & I &  0.5020 & 0.4657 & 0.4067 & 0.4286 &  0.4260 & 0.4320 & 0.4598 & 0.3983 & 0.4263 & 0.4243   \\
				& II & 0.9011 & 0.8495 & 0.8901 & 0.8952 & 0.8936 & 0.8847 & 0.8470 & 0.9185 & 0.8997 & 0.8972  \\ 
				& III & 0.9897 & 0.9967 & 0.9971 & 0.9969 &  0.9972 & 0.9991 & 0.9990 & 0.9976 & 0.9975 &  0.9975  \\ 
				& IV  &   1.0000 & 1.0000 & 1.0000 & 1.0000 & 1.0000 & 1.0000 & 1.0000 & 1.0000 & 1.0000 &  1.0000 \\ \hline \hline
				\multirow{4}{*}{\shortstack{after 960\\ time intants }} & I & 0.1404 & 0.1484 & 0.1483 & 0.1470 &  0.1489 &  0.1304 & 0.1375 & 0.1319 & 0.1403 &  0.1424 \\ 
				&	II & 0.4113 & 0.3817 & 0.3853 & 0.4231  & 0.4180 & 0.2997 & 0.3983 & 0.3875 & 0.4150 & 0.3995  \\ 
				&	III&   0.7108 & 0.7159 & 0.7726 & 0.7515 &  0.7471 &  0.7588 & 0.7959 & 0.7590 & 0.7525 &  0.7741 \\ 
				&	IV & 0.9914 & 0.9978 & 0.9933 & 0.9965 &  0.9952 &  0.9948 & 0.9953 & 0.9988 & 0.9985 & 0.9981 \\ \hline

			\end{tabular}
		\end{adjustbox}
		\caption{Simulated power of the pooled test for  a change in the standard deviations of the innovation processes with learning sample sizes $L_1=500$ and $L_2=1\,000$.}
		\label{pooledtestunknowncovpowercpinsd}
	\end{center}
\end{table}

The results of the second set of simulations (change in the AR coefficients) are provided in Tables \ref{squaredtestunknowncovpowercpincoefficients} and \ref{pooledtestunknowncovpowercpincoefficients}. Here the pooled test is preferable for early changes, as it has higher power than the SSQ test in this case. For a change in the middle of the sample or a late change the results are relatively comparable but in favor of the SSQ test when the learning sample is small. But for larger learning samples the pooled test is better for small dimensions.

\begin{table}
	\begin{center}
		\begin{adjustbox}{max width=\textwidth}
			\begin{tabular}{|c|c|c|c|c|c|c||c|c|c|c|c|c|c|c|c|} \hline
				&		& \multicolumn{5}{|c|}{$L_1 = 500$} & \multicolumn{5}{|c|}{$L_2 = 1\,000$} \\ \hline \hline
				Change	&	Case	& d=10 & d=50 & d=200 & d=500  & d=1000 & d=10 & d=50 & d=200 & d=500  & d=1000 \\ \hline 
				\multirow{4}{*}{\shortstack{after 240 \\ time instants}} & I &  0.1262 & 0.1933 & 0.1008 & 0.1216 &  0.1231 & 0.1299 & 0.1231 & 0.1093 & 0.1147 &  0.1064    \\
				& II & 0.3912 & 0.3669 & 0.2600  &0.2874 &  0.3085 & 0.1607 & 0.4495 & 0.2639 & 0.2563 & 0.2556   \\ 
				& III &  0.9908 & 0.6046 & 0.5883 & 0.5543 &  0.6312 & 0.9572 & 0.7524 & 0.6088 & 0.6388 &  0.6469  \\ 
				&	IV &  0.9213 & 0.9299 & 0.8627 & 0.9273 &  0.9402 & 0.9364 & 0.9894 & 0.9633 & 0.9061 &  0.9277   \\ \hline \hline
				\multirow{4}{*}{\shortstack{after 600 \\ time instants}} & I & 0.9666 & 0.9321 & 0.7871 & 0.8093 &  0.8190    & 0.5886 & 0.7900 & 0.8149 & 0.8184 &  0.8064  \\
				& II &  0.9999 & 0.9950 & 0.9993 & 0.9969 &  0.9981   &  0.9996 & 0.9991 & 0.9986 & 0.9973 &  0.9970  \\ 
				& III & 1.0000 & 1.0000 & 1.0000 & 1.0000  & 1.0000   & 1.0000 & 1.0000 & 1.0000 & 1.0000 &  1.0000 \\ 
				& IV  & 1.0000 & 1.0000 & 1.0000 & 1.0000 & 1.0000  &  1.0000 & 1.0000 & 1.0000 & 1.0000 &  1.0000 \\ \hline \hline
				\multirow{4}{*}{\shortstack{after 960\\ time intants }} & I & 0.9977 & 0.9985 & 0.9947 & 0.9941 &  0.9953 & 0.9998 & 0.9893 & 0.9958 & 0.9951 &  0.9941 \\ 
				&	II & 1.0000 & 1.0000 & 1.0000 & 1.0000 &  1.0000 & 1.0000 & 1.0000 & 1.0000 & 1.0000 &  1.0000  \\ 
				&	III& 1.0000 & 1.0000 & 1.0000 & 1.0000  & 1.0000 &  1.0000 & 1.0000 & 1.0000 & 1.0000 & 1.0000  \\ 
				&	IV & 1.0000 & 1.0000 & 1.0000 & 1.0000 &  1.0000  & 1.0000 & 1.0000 & 1.0000 & 1.0000  & 1.0000  \\ \hline

			\end{tabular}
		\end{adjustbox}
		\caption{Simulated power of the sum of squared errors test for a change in the coefficients with learning sample sizes $L_1=500$ and $L_2=1\,000$.}
		\label{squaredtestunknowncovpowercpincoefficients}
	\end{center}
\end{table}

\begin{table}
	\begin{center}
		\begin{adjustbox}{max width=\textwidth}
			\begin{tabular}{|c|c|c|c|c|c|c||c|c|c|c|c|c|c|c|c|} \hline
				&		& \multicolumn{5}{|c|}{$L_1 = 500$} & \multicolumn{5}{|c|}{$L_2 = 1\,000$} \\ \hline \hline
				Change	&	Case	& d=10 & d=50 & d=200 & d=500  & d=1000 & d=10 & d=50 & d=200 & d=500  & d=1000 \\ \hline 
				\multirow{4}{*}{\shortstack{after 240 \\ time instants}} & I & 0.5368 & 0.3360 & 0.3901 & 0.4405 &  0.4285 & 0.6967 & 0.4582 & 0.3719 & 0.3771 &  0.4286  \\
				& II &  0.9181 & 0.9334  &0.8335 & 0.8146  & 0.7951 &  0.9098 & 0.7745 & 0.8017 & 0.8003  & 0.8303  \\ 
				& III & 0.9981 & 0.9947 & 0.9948 & 0.9808  & 0.9862  & 0.9738 & 0.9909 & 0.9789 & 0.9848 &  0.9927   \\ 
				&	IV &  1.0000 & 1.0000 & 1.0000 & 1.0000 & 1.0000 &  0.9960 & 1.0000 & 1.0000 & 1.0000  & 0.9999   \\ \hline \hline
				\multirow{4}{*}{\shortstack{after 600 \\ time instants}} & I &  0.8127 & 0.8344 & 0.7587 & 0.7909 &  0.7794 & 0.8926 & 0.8150 & 0.7393 & 0.8334 & 0.7951      \\
				& II &  1.0000 & 0.9999 & 0.9980 & 0.9962 &  0.9961 & 0.9895 & 0.9864 & 0.9990 & 0.9980 & 0.9985    \\ 
				& III &  0.9999 & 1.0000 & 1.0000 & 1.0000 &  1.0000 &  1.0000 & 1.0000 & 1.0000 & 1.0000 &  1.0000  \\ 
				& IV  &  1.0000 & 1.0000 & 1.0000 & 1.0000 & 1.0000  &  1.0000 & 1.0000 & 1.0000 & 1.0000 & 1.0000  \\ \hline \hline
				\multirow{4}{*}{\shortstack{after 960\\ time intants }} & I & 0.9292 & 0.9469 & 0.9421 & 0.9362 & 0.9556 &  0.9976 & 0.9881 & 0.9687 & 0.9625 &  0.9656   \\ 
				&	II &  0.9998 & 1.0000 & 1.0000 & 0.9999 & 1.0000 & 1.0000 & 1.0000 & 1.0000 & 1.0000 & 1.0000   \\ 
				&	III& 1.0000 & 1.0000 & 1.0000 & 1.0000 & 1.0000 & 1.0000 & 1.0000 & 1.0000  &1.0000 &  1.0000     \\ 
				&	IV &  1.0000 & 1.0000 & 1.0000 & 1.0000 &  1.0000 & 1.0000 & 1.0000 & 1.0000 & 1.0000 &  1.0000  \\ \hline
				
			\end{tabular}
		\end{adjustbox}
		\caption{Simulated power of the pooled test statistic for a change in the coefficients  with learning sample sizes $L_1=500$ and $L_2=1\,000$.}
		\label{pooledtestunknowncovpowercpincoefficients}
	\end{center}
\end{table}

\clearpage

\subsection*{In-Sample Long-Run Variance Estimation}

Lastly, we investigate the behavior of the tests when no learning sample is available. Then one may rely on in-sample estimation of the unknown long-run variance parameter ignoring a possible change-point. This may lead to inconsistent estimation under the alternative hypothesis, but is valid under the null hypothesis.

The simulations were performed under the same conditions as above, but this time we only present the results for a change occurring after $960$ time instants in the standard deviations of the innovation processes. 

The estimated error of first kind for the sum of squared error test and the pooled test are given in Tables \ref{insamplelrvsquaredtestalpha} and \ref{insamplelrvpooledtestalpha}. One can see that now both tests keep the level very well.

\begin{table}[h]
	\begin{center}
		\begin{adjustbox}{max width=\textwidth}
			\begin{tabular}{|c|c|c|c|c|c|} \hline
				Case	& d=10 & d=50 & d=200 & d=500  & d=1000 \\ \hline \hline
				I   & 2.27 &  3.14 &   2.95 &  2.86 &   2.80 \\ \hline
				II  & 3.58 & 3.29 &  3.58 &  3.61  &  3.62 \\ \hline
				III & 4.21 & 3.83 &  4.02 &  4.14  &  4.00 \\ \hline
				IV  & 4.33 & 4.58 &  4.43  & 4.94  &  4.58 \\ \hline
			\end{tabular}
		\end{adjustbox}
		\caption{Estimated type I error rates for the sum of squared errors test for in-sample long-run variance parameter estimation}
		\label{insamplelrvsquaredtestalpha}
	\end{center}
\end{table}

\begin{table}[h]
	\begin{center}
		\begin{adjustbox}{max width=\textwidth}
			\begin{tabular}{|c|c|c|c|c|c|} \hline
				Case	& d=10 & d=50 & d=200 & d=500  & d=1000 \\ \hline \hline
				I & 3.10 & 3.38 & 3.77 & 3.64 & 3.58  \\ \hline
				II & 4.41 & 4.14 & 4.05 & 4.07 & 4.35 \\ \hline
				III & 4.88 &  4.72  &  4.42  &  4.28 &   4.42  \\ \hline
				IV & 4.71 & 4.45 & 4.75 &4.44  &4.62 \\ \hline
				
			\end{tabular}
		\end{adjustbox}
		\caption{Estimated type I error rates for the pooled test for in-sample long-run variance parameter estimation}
		\label{insamplelrvpooledtestalpha}
	\end{center}
\end{table}

The simulated power is given in Tables \ref{insamplelrvsquaredtestpower} and \ref{insamplelrvpooledtestpower}. It is visible that both tests have higher power compared to the tests using a learning sample. We may conclude that, for the cases investigated in this simulation study, in-sample estimation works very well and even leads to tests with higher detection power.

\begin{table}[h]
	\begin{center}
		\begin{adjustbox}{max width=\textwidth}
			\begin{tabular}{|c|c|c|c|c|c|} \hline
				Case	& d=10 & d=50 & d=200 & d=500  & d=1000 \\ \hline \hline
				I  & 0.6833 & 0.7347 & 0.7333 & 0.7351 & 0.7154 \\ \hline
				II & 0.9974 & 0.9962 & 0.9953 & 0.9960 & 0.9955 \\ \hline
				III & 1.0000 & 1.0000 & 0.9999 & 0.9999 &  1.0000 \\ \hline
				IV & 1.0000 & 1.0000 & 1.0000 & 1.0000 & 1.0000 \\ \hline
			\end{tabular}
		\end{adjustbox}
		\caption{Simulated power of the sum of squared errors test for a change in the standard deviations of the innovation processes and in-sample long-run variance parameter estimation}
		\label{insamplelrvsquaredtestpower}
	\end{center}
\end{table}

\clearpage

\begin{table}[h]
	\begin{center}
		\begin{adjustbox}{max width=\textwidth}
			\begin{tabular}{|c|c|c|c|c|c|} \hline
				Case	& d=10 & d=50 & d=200 & d=500  & d=1000 \\ \hline \hline
				I & 0.2093 & 0.2250 & 0.2321 & 0.2230 & 0.2314 \\ \hline
				II & 0.5083 & 0.4540 & 0.4356 & 0.4600 & 0.4487 \\ \hline
				III & 0.7779 & 0.7546 & 0.7603 & 0.7743&  0.7550  \\ \hline
				IV & 0.9999 & 0.9982 & 0.9992 & 0.9987 & 0.9986 \\ \hline
			\end{tabular}
		\end{adjustbox}
		\caption{Simulated power of the pooled test for a change in the standard deviations of the innovation processes and in-sample long-run variance parameter estimation}
		\label{insamplelrvpooledtestpower}
	\end{center}
\end{table}

\section{DATA ANALYSIS}

We applied the proposed methods to a  dataset on condition monitoring of a hydraulic systems publicly available  from the UCI Machine Learning Repository, see \cite{Helwig1}.

The data has been collected from an hydraulic test rig that consists of a main working station which is connected to a cooling circuit through an oil tank. Within the main working station oil is pumped through a system of accumulators, valves and containers. Various sensors measure  process quantities such as the pressure, volume flow or the temperature. Here we are only interested in the data of the pressure sensors of which there are $K=6$ in total, three within the main working station and three within the cooling circuit.

The test rig was constructed to simulate reversible fault conditions of various system components and the data corresponds to an experiment where such a fault condition was present. Clearly, a simulated fault in the accumulators or the valves will have an effect on the pressure within the system. We therefore analyzed the pressure sensor measurements to illustrate the proposed tests. In total  $\tilde{N} = 2205$ load cycles were run, each one lasting a minute. The pressure sensors sample data at a frequency of $100$ Hz per second, such that within the one minute time interval each pressure sensor is generating a time series of length $d = 6\,000$.  The first $211$ load cycles were simulated under the same conditions, and afterwards component faults were added about every 10 load cycles. Therefore, we used the first $L=110$ data points from each pressure sensor as our learning sample.

\begin{figure}[!ht]
	\begin{center}
		\includegraphics[scale=0.5]{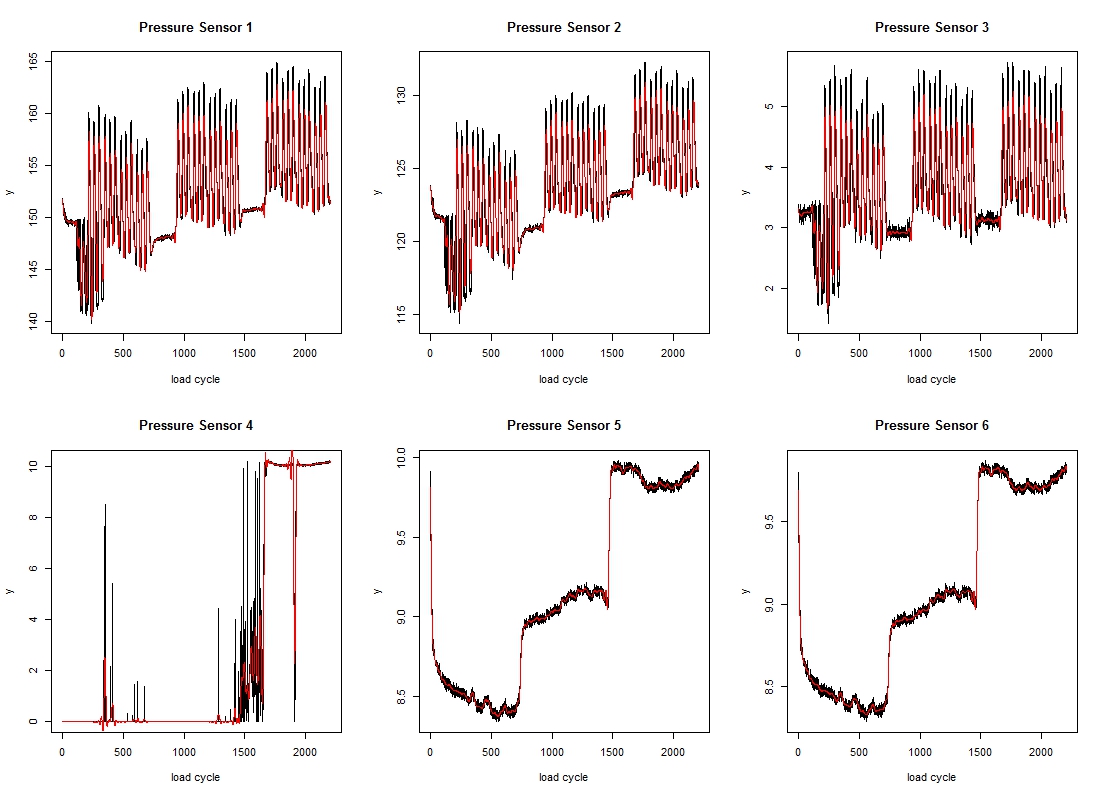} 
	\end{center}	
	\caption{Pressure sensor data (black) and cubic spline smoother (red) 10 seconds into each load cycle}
	\label{pressuredataandsplines}
	
\end{figure}

In a first preprocessing step residuals for the pressure sensor data were calculated. A cubic spline smoother for $f$, see \cite{GreenSilverman1993}, Sec. 2.3, was calculated and used to correct for the mean.
Figure \ref{pressuredataandsplines} shows the data time series at $\nu=1\,000$, which represents the pressure within the system exactly $10$ seconds into each load cycle. The black graph is the original data series, the red graph is the cubic spline smoother. We can see that the variance within the data generated by the pressure sensors within the main working station (pressure sensor 1 - 3) behaves similarly across all of the 2205 load cycles, but is completely different compared to the data generated by the pressure sensors within the cooling circuit.

To obtain $\ell_1 $-sparse projection vectors we calculated sparse principle directions following the approach of  \cite{Erichson2018} and given by the optimization problem
\begin{align*} \min \limits_{\matA, \matB} \quad g(A,B) = \frac12 \| \matS_L - \matS_L \matB \matA' \|_{\text{F}}^2 + \psi(\matB) , \qquad \text{s.t.} \qquad \matA' \matA = I,   \end{align*}
where $\matS_L \in \mathbb R^{K\cdot L \times d}$ is the sample variance-covariance matrix of the combined learning sample data, $\matA \in \mathbb R^{d \times r}$ is an orthonormal matrix, $\matB \in \mathbb R^{d \times r}$ is the sparse loadings matrix and $r \in \mathbb N$ is the number of sparse loading vectors to estimate (set to $r=5$). Further, an elastic net penalty,
$ \psi(\matB) =  \alpha \|\matB\|_1 + \frac12 \beta \|\matB\|_2^2 $
was used with $\alpha = \beta = 0.0001$ (the default values used by the \textsf{spca} function in the \textsf{R} package \textit{sparsepca}). Using the matrix $\matS_L$ allows us to use the same weighting vector for all of the $K=6$ pressure sensor datasets. 

Table \ref{sparseweightingvector} lists the first five leading sparse loading vectors $\vecw \in \mathbb R^d$ which explain, in total, $68.4\%$ of the variance within the learning sample data.	The vector $\vecw_5$ takes about half of the pressure sensor data generated within the one-minute load cycle into consideration, while for the others this is much less. It is clear that changes in the variance structure of projections based on these weighting vectors might not be detectable if too much data is omitted since fault conditions were simulated right at the start of the one-minute load cycles and sensors might not have recognized these changes within the time intervals where our vectors have nonzero entries.

\begin{table}[h]	
	\begin{center}
		\begin{tabular}{|c|c|c|c|c|} \hline
			Vector & Explained variance (in \%) & $\| \bullet\|_0$  \\ \hline
			$\vecw_1$ & 43.5  & 1080  \\ \hline
			$\vecw_2$ & 9.4  & 329  \\ \hline
			$\vecw_3$ & 7.4  & 340  \\ \hline
			$\vecw_4$ & 5.6  & 202  \\ \hline
			$\vecw_5$ & 2.5  & 3100  \\ \hline
		\end{tabular}
		\caption{Explained variance and $\ell_0$ - norm of the first five sparse loading vectors}
		\label{sparseweightingvector}
	\end{center}
\end{table}

\begin{figure}[!ht]
	\begin{center}
		\includegraphics[scale=0.6]{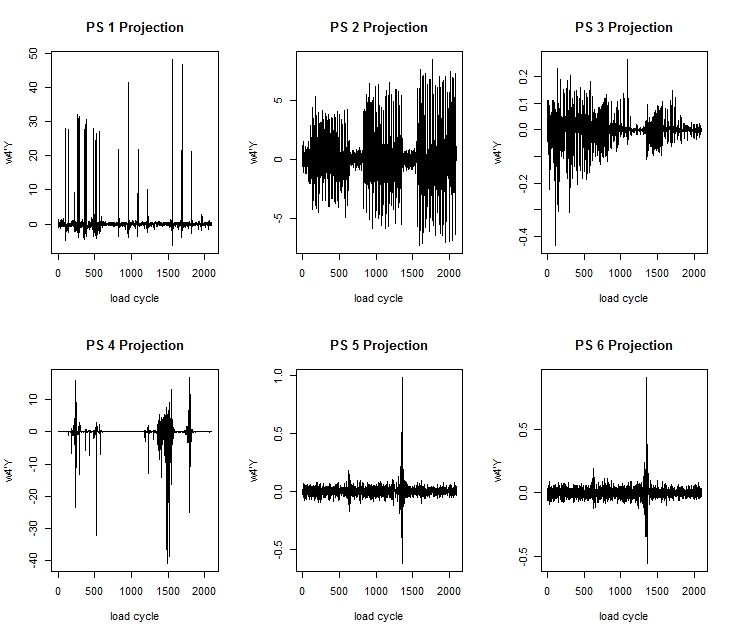} 
		\caption{Projections with weighting vector $\vecw_4$}
		\label{w4projections}
	\end{center}	
\end{figure}

Figure \ref{w4projections} shows the projections with the weighting vector $\vecw_4$, where changes in the variance structure of the data generated by the second pressure sensor data are clearly visible.

Finally, Table \ref{pressuresensorresults} shows the results of our change-point tests.
It turns out that using the sum of squared errors test with the test statistic $\breve{Q}$ we can always verify the existence of a change within the projected data set in all of the five cases, while the test based on the pooled sample variance-covariance matrix with test statistic $\breve{V}$ is only unable to detect a change in the case of $\vecw_4$. This is due to the fact that the majority of nonzero entries of $\vecw_4$ lie within the first 2000 components of the vector while there are no values taking the second half of each load cycle into consideration.

\begin{table}[h]
	\begin{center}
		\begin{tabular}{|c||c|c||c|c|c|c|} \hline
			Vector & $\breve{Q}$ & $c_{\text{crit}}^{(\breve{Q})}$ & $\breve{V}$ & $c_{\text{crit}}^{(\breve{V})}$ \\ \hline
			$\vecw_1$ & $8.24 \cdot 10^{22}$  & 7.08 & 311622.75 & 22638.27 \\ \hline
			$\vecw_2$ & $ 2.52 \cdot 10^{22} $  & 7.08 & 11147.93 &4465.77 \\ \hline
			$\vecw_3$ & $ 1.43 \cdot 10^{23} $ & 7.08 & 6430.08 & 6240.02 \\ \hline
			$\vecw_4$ & $ 5.86 \cdot 10^{20} $ & 7.08 & 95.31 & 5154.20 \\ \hline
			$\vecw_5$ & $ 2.14 \cdot 10^{22} $ & 7.08 & 19820.19 & 1524.12 \\ \hline
		\end{tabular}
		\caption{Change-point test results}
		\label{pressuresensorresults}
	\end{center}
\end{table}

\section{PROOFS} 

\textit{Proof of Theorem \ref{squaredtheorem}:}\\

Due to Theorem 4.1 in \cite{steland1}, for any $T>0$ and each $j=1,\dots,K$ there exists an equivalent version of $D_{T,k_j}^{(j)}, k_j \in \mathbb N,$ and thus of $\mathcal{D}_T^{(j)}(t_j)\,, t_j \geq 0,$ denoted by $\tilde{D}_{T,k_j}^{(j)}$ and $\tilde{\mathcal{D}}_T^{(j)}(t_j)$, and a standard Brownian motion $\{\tilde{B}_{T}^{(j)}(t_j)\, : \, t_j \geq 0 \}$ on a new probability space $(\tilde{\Omega}^{(j)},\tilde{\mathcal{F}}^{(j)},\tilde{P}^{(j)})$, such that
\begin{align} \label{DTjtapproxx} \left \vert \tilde{D}_{T,k_j}^{(j)} - \alpha_{T,j} \tilde{B}_{T}^{(j)}(k_j)       \right \vert \leq C_{T,j} k_j^{\frac{1}{2} - \lambda_j} , \quad t_j > 0, \, j=1,\dots,K,    \end{align}
$\tilde{P}^{(j)}$-a.s., where $C_{T,j} <  \infty$ are constants and $\lambda_j > 0$ for $j=1,\dots,K$. If additionally $C_{T,j} N_j^{-\lambda_j} = o(1)$ for $j \in \{1,\dots,K\}$ as $T \to \infty$, then we also have the strong approximation 
\begin{align} \label{cadlagapproxx} \sup \limits_{t_j \in [0,1]} \left \vert  \tilde{\mathcal{D}}_{T}^{(j)} (t_j) - \alpha_{T,j} \tilde{B}_{T}^{(j)} \left ( \frac{\lfloor t_j N_j \rfloor}{N_j}  \right )  \right \vert  = o(1), \quad \tilde{P}^{(j)}-\text{a.s.}, \quad T \to \infty,    \end{align}
for the $j$-th \cadlag version of the bilinear form based on the equivalent version of the centered and appropriately scaled sample variance-covariance matrix $\hat{\tilde{\bm{\Sigma}}}_{T}^{(j)}$. Similarly, for the bridge process $\tilde{\Delta}_T^{(j)}$ we have the strong approximation
\begin{align} \label{bridgeapproxx} \sup \limits_{t_j \in [0,1]} \left \vert \tilde{\Delta}_{T}^{(j)}(t_j) - \alpha_{T,j} \overline{\tilde{B}}_{T}^{(j)} \left ( \frac{ \lfloor t_j N_j \rfloor}{N_j}  \right )  \right \vert = o(1), \quad \tilde{P}^{(j)}-\text{a.s.}, \quad T \to \infty,  \end{align}
where $\overline{\tilde{B}}_{T}^{(j)}(t_j) = \tilde{B}_{T}^{(j)}(t_j) - t_j \tilde{B}_{T}^{(j)}(1), t_j \in [0,1], j=1,\dots,K$, are Brownian bridges. 

It is important to note that this construction can be done in such a way that, under the measure $ \tilde{P} $ on the cartesian product $ \times_{j=1}^K \tilde{\Omega}_j $, the coordinates are independent. Heuristically, this is clear, but the details are subtle, since iteratively applying the construction for each coordinate does not guarantee independence of the Brownian motions. Here, a second step is needed. Indeed, the proof of \cite{philipp}, which refers to \cite{PhilippStout1975} and \cite{BerkesPhilipp1979} works as follows: First, one embeds the summands of a partial sum, $ S_n $, into a Brownian motion, $B$, defined on $ ([0,1], \lambda) $, where $ \lambda $ denotes Lebesgue measure, and redefines the triple consisting of the partial sum, its approximation, $S_n'$, and the Brownian motion on a further probability space, such that the joint distributions $ (S_n', B) $ and $ (S_n, S_n') $ remain unchanged, see \cite[Lemma~A.1, proof of Th.~2]{BerkesPhilipp1979} and \cite[p.23]{PhilippStout1975}. Start with  $ \bar{\Omega} = \times_{j=1}^K \bar{\Omega}^{(j)} $ with $ \bar{\Omega}^{(j)} := [0,1] $ and equip $ \bar{\Omega} $ with the product measure $ \bar{P} = \times_{j=1}^K \lambda^{(j)} $ of $ K$ Lebesgue measures $ \lambda^{(j)} = \lambda $ on $ [0,1] $. Now take the $j$ th pair $ ([0,1], \lambda_j) $ for the above construction. The resulting $K$ triples can now be constructed on a new product probability space, such that the triples are independent and each triple is as required. This gives $ \tilde{\Omega} = \times_{j=1}^K \tilde{\Omega}^{(j)} $ on which the equivalent processes and the $K$ independent Brownian motions are defined.

Note that the $j$-th process $\tilde{\mathcal{D}}_T^{(j)}$ and the $j$-th Brownian motion $\tilde{B}_T^{(j)}$ only depend on the $j$-th parameter $t_j$ and $\lfloor t_j N_j \rfloor \slash N_j $ respectively and not on any other other parameter $t_{\tilde{j}}$ for $\tilde{j} \in \{1,\dots,K\}$ with $\tilde{j} \neq j$, which allows us to express the supremum $\sup_{\vect \in [0,1]^K}$ of sums running over $j=1,\dots,K$ as sums of the $K$ suprema $\sup_{t_j \in [0,1]}$.
\vspace{0.2cm}  
\\Now, using the inequality $ \vert a^2 - b^2 \vert  \leq \vert a-b \vert^2 + 2 \vert a-b\vert \, \vert b \vert $ for $a,b \in \mathbb R$ we have for $k_j \in \mathbb N$
\begin{align*}  \left \vert  \left ( \tilde{D}_{T,k_j}^{(j)} \right )^2 - \left ( \alpha_{T,j} \tilde{B}_T^{(j)} (k_j) \right )^2  \right \vert & \leq  \left \vert \tilde{D}_{T,k_j}^{(j)} - \alpha_{T,j} \tilde{B}_{T}^{(j)}(k_j)       \right \vert^2 + 2  \left \vert \tilde{D}_{T,k_j}^{(j)} - \alpha_{T,j} \tilde{B}_{T}^{(j)}(k_j)       \right \vert \, \left \vert  \alpha_{T,j} \tilde{B}_T^{(j)} (k_j) \right \vert.  \end{align*}
Recalling the absolute deviation $ E | \tilde{B}_T^{(j)}(k_j) | = \sqrt{\frac{2 k_j}{\pi}} $, such that
$ \vert \alpha_{T,j} \tilde{B}_T^{(j)}(k_j) \vert = \mathcal{O}_P \left ( \alpha_{T,j} k_j^{\frac12} \right ) $, we obtain for all $k_j \in \mathbb N$ and $j=1,\dots,K$,
\begin{align} \nonumber  \left \vert  \left ( \tilde{D}_{T,k_j}^{(j)} \right )^2 - \left ( \alpha_{T,j} \tilde{B}_T^{(j)} (k_j) \right )^2  \right \vert \, \, \, & \stackrel{\hidewidth \eqref{DTjtapproxx} \hidewidth}{=} \, \, \, \mathcal{O} \left (C_{T,j}^2 k_j^{1-2\lambda_j} \right ) + \mathcal{O} \left (C_{T,j} k_j^{\frac12 - \lambda_j} \right ) \mathcal{O}_P \left (\alpha_{T,j} k_j^{\frac12} \right ) \\ & \stackrel{\hidewidth \hidewidth}{=} \, \, \, \mathcal{O}_P \left (\alpha_{T,j} \tilde{C}_{T,j} k_j^{1 - \lambda_j} \right ) \label{squaredbound} , \end{align}
$\tilde{P}^{(j)}-$ a.s., with $\tilde{C}_{T,j} = \max \{ C_{T,j}^2, C_{T,j}\}$ and therefore by the triangle inequality
\begin{align*} \left \vert  \sum_{j=1}^K \left ( \tilde{D}_{T,k_j}^{(j)}  \right )^2 - \sum_{j=1}^K \left ( \alpha_{T,j} \tilde{B}_T^{(j)}(k_j) \right )^2 \right \vert &% \leq \sum_{j=1}^K \left \vert \left ( D_{T,k_j}^{(j)}  \right )^2 - \left ( \alpha_{T,j} B_T^{(j)}(k_j) \right )^2 \right \vert \\ &
\stackrel{\hidewidth \eqref{squaredbound} \hidewidth}{=} \, \, \, \sum_{j=1}^K \mathcal{O}_P \left (\alpha_{T,j} \tilde{C}_{T,j} k_j^{1 - \lambda_j} \right )   \end{align*}
$\tilde{P}-$ a.s., follows. Now, since
\begin{align*} \tilde{\mathcal{D}}_T^{(j)}(t_j) = \frac{1}{\sqrt{N_j}} \vecv_{T,j}' \left ( \hat{\tilde{\bm{\Sigma}}}_{T,\lfloor t_j N_j \rfloor}^{(j)} - \bm{\Sigma}_{T,\lfloor t_j n_j \rfloor}^{(j)} \right ) \vecw_{T,j} = \frac{1}{\sqrt{N_j}} \tilde{D}_{T,\lfloor t_j N_j \rfloor}^{(j)}, \quad t_j \geq 0, \end{align*}
for all $j=1,\dots,K$, we have due to the scaling property of Brownian motion
\begin{align*} & \sup \limits_{\vect \in [0,1]^K} \left \vert \sum_{j=1}^K \left ( \tilde{\mathcal{D}}_T^{(j)}(t_j)  \right )^2 - \sum_{j=1}^K \left ( \alpha_{T,j} \tilde{B}_T^{(j)} \left ( \frac{\lfloor t_j N_j \rfloor}{N_j} \right ) \right )^2  \right \vert \\ & \leq \sum_{j=1}^K \sup \limits_{t_j \in [0,1]} \left \vert \left (  \tilde{\mathcal{D}}_T^{(j)}(t_j)   \right )^2 - \left ( \alpha_{T,j} \tilde{B}_T^{(j)} \left ( \frac{\lfloor t_j N_j \rfloor}{N_j} \right )  \right )^2     \right \vert \\ & = \sum_{j=1}^K \sup \limits_{t_j \in [0,1]} \left \vert \left (  \frac{1}{\sqrt{N_j}} \tilde{D}_{T,\lfloor t_j N_j \rfloor}^{(j)}   \right )^2 - \left ( \alpha_{T,j} \tilde{B}_T^{(j)} \left ( \frac{\lfloor t_j N_j \rfloor}{N_j} \right )  \right )^2     \right \vert \\ & =  \sum_{j=1}^K \frac{1}{N_j}  \sup \limits_{t_j \in [0,1]} \left \vert \left (  \tilde{D}_{T,\lfloor t_j N_j \rfloor}^{(j)}   \right )^2 - \left ( \alpha_{T,j} \sqrt{N_j} \tilde{B}_T^{(j)} \left ( \frac{\lfloor t_j N_j \rfloor}{N_j} \right )  \right )^2     \right \vert  \\ & \stackrel{d}{=}  \sum_{j=1}^K \frac{1}{N_j}  \sup \limits_{t_j \in [0,1]} \left \vert \left (  \tilde{D}_{T,\lfloor t_j N_j \rfloor}^{(j)}   \right )^2 - \left ( \alpha_{T,j}  \tilde{B}_T^{(j)} \left (\lfloor t_j N_j \rfloor \right )  \right )^2     \right \vert. \end{align*}
Using \eqref{squaredbound} we then have
\begin{align*}  & \sup \limits_{\vect \in [0,1]^K} \left \vert \sum_{j=1}^K \left ( \tilde{\mathcal{D}}_T^{(j)}(t_j)  \right )^2 - \sum_{j=1}^K \left ( \alpha_{T,j} \tilde{B}_T^{(j)} \left ( \frac{\lfloor t_j N_j \rfloor}{N_j} \right ) \right )^2  \right \vert  = \sum_{j=1}^K \frac{1}{N_j} \sup \limits_{t_j \in [0,1]} \mathcal{O}_P \left ( \alpha_{T,j} \tilde{C}_{T,j} \lfloor t_j N_j \rfloor^{1 - \lambda_j} \right ) \\ & = \sum_{j=1}^K \frac{1}{N_j} \alpha_{T,j} \mathcal{O}_P \left (  \tilde{C}_{T,j} N_j^{1 - \lambda_j} \right ) = \sum_{j=1}^K \alpha_{T,j} \mathcal{O}_P \left ( \tilde{C}_{T,j} N_j^{- \lambda_j} \right ),   \end{align*}
such that \eqref{squaredapprox1} follows when additionally $\alpha_{T,j} \to \alpha_j^*$ and $\tilde{C}_{T,j}N_j^{ - \lambda_j} = o(1)$ for all $j=1,\dots,K$ as $T \to \infty$. Further, 
\begin{align*}  &  \sup \limits_{\vect \in [0,1]^K} \left \vert \sum_{j=1}^K \left ( \tilde{\Delta}_T^{(j)}(t_j) \right )^2 - \sum_{j=1}^K \left [ \alpha_{T,j} \tilde{\overline{B}}_T^{(j)} \left ( \frac{\lfloor t_j N_j \rfloor}{N_j} \right )     \right ]^2  \right \vert \\ & \leq \sup \limits_{\vect \in [0,1]^K} \left \vert \sum_{j=1}^K  \left ( \tilde{\mathcal{D}}_T^{(j)} \left ( \frac{\lfloor t_j N_j \rfloor}{N_j} \right )   \right )^2 - \sum_{j=1}^K  \left ( \alpha_{T,j} \tilde{B}_T^{(j)} \left ( \frac{\lfloor t_j N_j \rfloor}{N_j} \right )  \right )^2    \right \vert \\ & \quad + 2 \sup \limits_{\vect \in [0,1]^K} \left \vert \sum_{j=1}^K   \tilde{\mathcal{D}}_T^{(j)} \left ( \frac{\lfloor t_j N_j \rfloor}{N_j} \right ) \frac{\lfloor t_j N_j \rfloor}{N_j} \tilde{\mathcal{D}}_T^{(j)}(1) - \sum_{j=1}^K \alpha_{T,j}^2 \tilde{B}_T^{(j)} \left ( \frac{\lfloor t_j N_j \rfloor}{N_j} \right ) \frac{\lfloor t_j N_j \rfloor}{N_j} \tilde{B}_T^{(j)}(1) \right \vert \\ & \quad + \sup \limits_{\vect \in [0,1]^K} \left \vert \sum_{j=1}^K  \left (  \frac{\lfloor t_j N_j \rfloor}{N_j} \tilde{\mathcal{D}}_T^{(j)}(1)    \right )^2  - \sum_{j=1}^K   \left ( \alpha_{T,j} \frac{\lfloor t_j N_j \rfloor}{N_j} \tilde{B}_T^{(j)}(1)  \right )^2  \right \vert \\ & = I_{T,1} + I_{T,2} + I_{T,3}. \end{align*}
If $\alpha_{T,j} \to \alpha_j^*$ and $\tilde{C}_{T,j} N_j^{-\lambda_j} = o(1)$ as $T \to \infty$ for $j=1,\dots,K,$  we have $I_{T,1} = o_P(1)$ by \eqref{squaredapprox1} and 
\begin{align*} I_{T,3} & \leq \sum_{j=1}^K \sup \limits_{t_j \in [0,1]}  \left ( \frac{\lfloor t_j N_j \rfloor}{N_j}   \right )^2  \left \vert \left ( \tilde{\mathcal{D}}_T^{(j)}(1)  \right )^2 - \left ( \alpha_{T,j} \tilde{B}_T^{(j)}(1) \right )^2    \right \vert \\ & =  \sum_{j=1}^K \sup \limits_{t_j \in [0,1]}  \left ( \frac{\lfloor t_j N_j \rfloor}{N_j}   \right )^2 \frac{1}{N_j} \left \vert  (\tilde{D}_{T,N_j}^{(j)})^2 - (\alpha_{T,j} \sqrt{N_j} \tilde{B}_T^{(j)}(1))^2 \right \vert \\ & \stackrel{d}{=} \sum_{j=1}^K \frac{1}{N_j} \left \vert  (\tilde{D}_{T,N_j}^{(j)})^2 - (\alpha_{T,j} \tilde{B}_T^{(j)}(N_j))^2 \right \vert \\ & \stackrel{\hidewidth \eqref{squaredbound} \hidewidth}{=}  \, \, \, \sum_{j=1}^K \frac{1}{N_j} \alpha_{T,j} \mathcal{O}_P \left (\tilde{C}_{T,j} N_j^{1 - \lambda_j} \right ) = o_P(1), \quad T \to \infty.  \end{align*}
Lastly, using the inequality  $\vert ab - cd \vert \leq \vert a-c \vert \, (\vert b -d \vert + \vert d \vert ) + \vert c \vert \, \vert b-d \vert$ for $a,b,c,d \in \mathbb R$ we get
\begin{align*}  I_{T,2} &=  2 \sup \limits_{\vect \in [0,1]^K} \left \vert \sum_{j=1}^K   \tilde{\mathcal{D}}_T^{(j)} \left ( \frac{\lfloor t_j N_j \rfloor}{N_j} \right ) \frac{\lfloor t_j N_j \rfloor}{N_j} \tilde{\mathcal{D}}_T^{(j)}(1) - \sum_{j=1}^K \alpha_{T,j}^2 \tilde{B}_T^{(j)} \left ( \frac{\lfloor t_j N_j \rfloor}{N_j} \right ) \frac{\lfloor t_j N_j \rfloor}{N_j} \tilde{B}_T^{(j)}(1) \right \vert \\ & \leq 2 \sum_{j=1}^K \sup \limits_{t_j \in [0,1]} \left \{ \left \vert  \tilde{\mathcal{D}}_T^{(j)} \left ( \frac{\lfloor t_j N_j \rfloor}{N_j} \right ) - \alpha_{T,j} \tilde{B}_T^{(j)} \left ( \frac{\lfloor t_j N_j \rfloor}{N_j} \right )      \right \vert \right. \\ & \qquad \qquad \qquad \times \left.  \left ( \left \vert  \frac{\lfloor t_j N_j \rfloor}{N_j} \tilde{\mathcal{D}}_T^{(j)}(1) - \alpha_{T,j} \frac{\lfloor t_j N_j \rfloor}{N_j} \tilde{B}_T^{(j)}(1)    \right \vert + \left \vert \alpha_{T,j} \frac{\lfloor t_j N_j \rfloor}{N_j} \tilde{B}_T^{(j)}(1)    \right \vert  \right ) \right \} \\ & \quad + 2 \sum_{j=1}^K \sup \limits_{t_j \in [0,1]} \left \vert \alpha_{T,j} \tilde{B}_T^{(j)} \left ( \frac{\lfloor t_j N_j \rfloor}{N_j} \right )  \right \vert \, \left \vert \frac{\lfloor t_j N_j \rfloor}{N_j} \tilde{\mathcal{D}}_T^{(j)}(1) -  \alpha_{T,j} \frac{\lfloor t_j N_j \rfloor}{N_j} \tilde{B}_T^{(j)}(1)   \right \vert \\ & \leq 2 \sum_{j=1}^K \sup \limits_{t_j \in [0,1]} \left \{ \left \vert  \tilde{\mathcal{D}}_T^{(j)} \left ( \frac{\lfloor t_j N_j \rfloor}{N_j} \right ) - \alpha_{T,j} \tilde{B}_T^{(j)} \left ( \frac{\lfloor t_j N_j \rfloor}{N_j} \right )      \right \vert \right. \\ & \qquad \qquad \qquad \times \left.  \left ( \sup \limits_{t_j \in [0,1]} \left \vert  \frac{\lfloor t_j N_j \rfloor}{N_j} \tilde{\mathcal{D}}_T^{(j)}(1) - \alpha_{T,j} \frac{\lfloor t_j N_j \rfloor}{N_j} \tilde{B}_T^{(j)}(1)    \right \vert + \sup \limits_{t_j \in [0,1]} \left \vert \alpha_{T,j} \frac{\lfloor t_j N_j \rfloor}{N_j} \tilde{B}_T^{(j)}(1)    \right \vert  \right ) \right \} \\ & + 2 \sum_{j=1}^K \sup \limits_{t_j \in [0,1]} \left \vert \alpha_{T,j} \tilde{B}_T^{(j)} \left ( \frac{\lfloor t_j N_j \rfloor}{N_j} \right )  \right \vert \, \sup \limits_{t_j \in [0,1]} \left \vert \frac{\lfloor t_j N_j \rfloor}{N_j} \tilde{\mathcal{D}}_T^{(j)}(1) -  \alpha_{T,j} \frac{\lfloor t_j N_j \rfloor}{N_j} \tilde{B}_T^{(j)}(1)   \right \vert. \end{align*}
Note that, since $ \tilde{\mathcal{D}}_T^{(j)}\left ( \frac{\lfloor t_j N_j \rfloor}{N_j}  \right ) =   \tilde{\mathcal{D}}_T^{(j)}(t_j)$ for all $t_j \in [0,1], j=1,\dots,K$, 
\begin{align*} & \sup \limits_{t_j \in [0,1]}  \left \vert  \tilde{\mathcal{D}}_T^{(j)} \left ( \frac{\lfloor t_j N_j \rfloor}{N_j} \right ) - \alpha_{T,j} \tilde{B}_T^{(j)} \left ( \frac{\lfloor t_j N_j \rfloor}{N_j} \right )      \right \vert = \sup \limits_{t_j \in [0,1]} \left \vert  \tilde{\mathcal{D}}_T^{(j)} \left ( t_j \right ) - \alpha_{T,j} \tilde{B}_T^{(j)} \left ( \frac{\lfloor t_j N_j \rfloor}{N_j} \right )      \right \vert \\ & \stackrel{\hidewidth \eqref{cadlagapproxx} \hidewidth}{=} \, \, \, \mathcal{O} \left (  C_{T,j} N_j^{-\lambda_j}  \right ). \end{align*}
Further, observing that $ \sup_{t_j \in [0,1]} \frac{\lfloor t_j N_j \rfloor}{N_j} = 1 $, we have
\begin{align*} & \sup \limits_{t_j \in [0,1]} \left \vert \frac{\lfloor t_j N_j \rfloor}{N_j} \tilde{\mathcal{D}}_T^{(j)}(1) -  \alpha_{T,j} \frac{\lfloor t_j N_j \rfloor}{N_j} \tilde{B}_T^{(j)}(1)   \right \vert = \left \vert \tilde{\mathcal{D}}_T^{(j)}(1) -  \alpha_{T,j}  \tilde{B}_T^{(j)}(1)   \right \vert  \stackrel{\eqref{cadlagapproxx}}{=} \mathcal{O} \left (C_{T,j}N_j^{-\lambda_j} \right ) \end{align*}
and $\sup_{t_j \in [0,1]} \left \vert  \alpha_{T,j} \frac{\lfloor t_j N_j \rfloor}{N_j} \tilde{B}_T^{(j)}(1)   \right \vert  = \mathcal{O}_P(\alpha_{T,j})$ as well as $\sup_{t_j \in [0,1]} \left \vert \alpha_{T,j} \tilde{B}_T^{(j)} \left ( \frac{\lfloor t_j N_j \rfloor}{N_j} \right )  \right \vert =\mathcal{O}_P(\alpha_{T,j})$ by Markov's inequality. Combining these facts we obtain
\begin{align*} I_{T,2} &= \sum_{j=1}^K \left ( \mathcal{O}\left (C_{T,j} N_j^{-\lambda_j} \right ) \left ( \mathcal{O}(C_{T,j} N_j^{-\lambda_j}) + \mathcal{O}_P(\alpha_{T,j})    \right ) + \mathcal{O}_P(\alpha_{T,j}) \mathcal{O} \left ( C_{T,j} N_j^{-\lambda_j}  \right )  \right ) \\ & = \sum_{j=1}^K \alpha_{T,j} \mathcal{O}(\tilde{C}_{T,j} N_j^{-\lambda_j}) \end{align*}
and finally arrive at 
\begin{align*}  &  \sup \limits_{\vect \in [0,1]^K} \left \vert \sum_{j=1}^K \left ( \tilde{\Delta}_T^{(j)}(t_j) \right )^2 - \sum_{j=1}^K \left [ \alpha_{T,j} \overline{\tilde{B}}_T^{(j)} \left ( \frac{\lfloor t_j N_j \rfloor}{N_j} \right )     \right ]^2  \right \vert = o_P(1), \quad T \to \infty,    \end{align*}
which concludes the proof.$ \hfill \Box$\\
\vspace{0.2cm}
\\\textit{Proof of Corollary \ref{squaredcorollary1}:}\\
\vspace{0.1cm}
\\When standardizing with $ \alpha_{T,j} $ the proof is straightforward and therefore omitted, see \cite{Mause2019}. Especially, one gets
\begin{align*}
\left| \left( \tilde{\mathcal{D}}_{T,j}(t_j) / \alpha_{T,j}  \right)^2 - \left(  \tilde{B}_T^{(j)}( \lfloor t_j N_j \rfloor / N_j ) \right)^2 \right| = \mathcal{O}_P( N_j^{-\lambda_j} ) \end{align*}
Now replace $ \alpha_{T,j} $ by  estimators with $ | \widehat{\alpha}_{T,j} / \alpha_{T,j} - 1 | = o_P(1) $. For brevity fix $j$ and denote $ \tilde{\mathcal{D}}_j = \tilde{\mathcal{D}}_{T,j}(t_j) $, $ \tilde{B}_j = \tilde{B}_T^{(j)}( \lfloor t_j N_j \rfloor / N_j )  $, $ \alpha_j = \alpha_{T,j} $ and $ \widehat{\alpha}_j = \widehat{\alpha}_{T,j} $. The proof can now be completed using 
\begin{align*}
\left| \frac{\tilde{\mathcal{D}}_j }{ \widehat{\alpha}_j } - \tilde{B}_j \right|
\le \left|  \frac{\alpha_j}{ \widehat{\alpha}_j } \right|
\left(  \left| \frac{\tilde{\mathcal{D}}_j}{\alpha_j} - \tilde{B}_j \right| + |\tilde{B}_j|\left| 1 - \frac{\widehat{\alpha}_j}{\alpha_j} \right| \right) = o_P(1), \end{align*}
as $ T \to \infty $, and the inequality 
\begin{align*}
\left| \left( \frac{\tilde{\mathcal{D}}_j }{ \widehat{\alpha}_j } \right)^2 - \tilde{B}_j^2 \right|
\le \left( \left| \frac{\tilde{\mathcal{D}}_j }{ \widehat{\alpha}_j } - \tilde{B}_j \right| + 2 | \tilde{B}_j |  \right) \left| \frac{\tilde{\mathcal{D}}_j }{ \widehat{\alpha}_j } - \tilde{B}_j \right|, \end{align*}
recalling the fact that, almost surely, $ \sup_{t_j \in [0,1]} | \tilde{B}_j(t_j) | < \infty $. 
$ \hfill \Box$\\
\vspace{0.2cm}
\\\textit{Proof of Corollary \ref{squaredcusumcorollary}:}\\
%\vspace{0.1cm}
\\Let $M \subset \mathbb R^p, \, p \in \mathbb N,$ be an arbitrary set and $g_1,g_2: M \mapsto \mathbb R$ be arbitrary functions. Then
\begin{align} 
\label{distancemaxbound} \left \vert  \sup_{x \in M} \vert g_1(x) \vert - \sup \limits_{x \in M} \vert g_2(x) \vert     \right \vert \leq \sup \limits_{x \in M} \vert g_1(x) - g_2(x)   \vert. 
\end{align}
Further,  $\sup_{t_j \in [0,1]} \left ( \frac{1}{\alpha_{T,j}} \tilde{\mathcal{D}}_T^{(j)}(t_j) \right )^2
= \max_{k_j \leq N_j} \left ( \frac{1}{\alpha_{T,j}} \tilde{\mathcal{D}}_T^{(j)} \left ( \frac{k_j}{N_j} \right )  \right )^2 $ as well as $  \sup_{t_j \in [0,1]} \left ( \tilde{B}_T^{(j)} \left ( \frac{\lfloor t_j N_j \rfloor}{N_j}  \right )  \right )^2 = \max_{k_j \leq N_j}  \left ( \tilde{B}_T^{(j)} \left ( \frac{k_j}{N_j}  \right )  \right )^2.$
We now obtain
\begin{align*}  & \left \vert \sum_{j=1}^K \max \limits_{k_j \leq N_j} \left ( \frac{1}{\alpha_{T,j}} \tilde{\mathcal{D}}_T^{(j)} \left ( \frac{k_j}{N_j}  \right )  \right )^2  - \sum_{j=1}^K \max \limits_{k_j \leq N_j} \left ( \tilde{B}_T^{(j)} \left ( \frac{k_j}{N_j} \right ) \right )^2 \right \vert  \\ & = \left \vert  \sum_{j=1}^K \sup \limits_{t_j \in [0,1]} \left ( \frac{1}{\alpha_{T,j}} \tilde{\mathcal{D}}_T^{(j)} \left ( t_j  \right )  \right )^2  - \sum_{j=1}^K \sup \limits_{t_j \in [0,1]} \left ( \tilde{B}_T^{(j)} \left ( \frac{\lfloor t_j N_j \rfloor}{N_j} \right ) \right )^2  \right \vert \\ & = \left \vert \sup \limits_{\vect \in [0,1]^K}  \sum_{j=1}^K  \left ( \frac{1}{\alpha_{T,j}} \tilde{\mathcal{D}}_T^{(j)} \left ( t_j  \right )  \right )^2  - \sup \limits_{\vect \in [0,1]^K} \sum_{j=1}^K \left ( \tilde{B}_T^{(j)} \left ( \frac{\lfloor t_j N_j \rfloor}{N_j} \right ) \right )^2  \right \vert \\ & \stackrel{\hidewidth \eqref{distancemaxbound} \hidewidth}{\leq} \sup \limits_{\vect \in [0,1]^K}  \left \vert   \sum_{j=1}^K  \left ( \frac{1}{\alpha_{T,j}} \tilde{\mathcal{D}}_T^{(j)} \left ( t_j  \right )  \right )^2  - \sum_{j=1}^K \left ( \tilde{B}_T^{(j)} \left ( \frac{\lfloor t_j N_j \rfloor}{N_j} \right ) \right )^2  \right \vert \\ & \stackrel{\hidewidth \eqref{squaredapproxst1} \hidewidth}{=} \, \, \, o_P(1)  \end{align*}
when $\tilde{C}_{T,j}N_j^{ - \lambda_j} = o(1)$ and $\alpha_{T,j} \to \alpha_j^* \in (0,\infty)$ for all $j=1,\dots,K$ as $T \to \infty$. \\
\vspace{0.1cm}
\\Similarly, 
\begin{align*} & \left \vert \sum_{j=1}^K \max \limits_{k_j \leq N_j} \left ( \frac{1}{\alpha_{T,j}} \tilde{\Delta}_T^{(j)} \left ( \frac{k_j}{N_j}  \right )  \right )^2  - \sum_{j=1}^K \max \limits_{k_j \leq N_j} \left ( \overline{\tilde{B}}_T^{(j)} \left ( \frac{k_j}{N_j} \right ) \right )^2 \right \vert \\ & = \left \vert \sup \limits_{\vect \in [0,1]^K} \sum_{j=1}^K  \left ( \frac{1}{\alpha_{T,j}} \tilde{\Delta}_T^{(j)} \left ( t_j  \right )  \right )^2  - \sup \limits_{\vect \in [0,1]^K} \sum_{j=1}^K  \left ( \overline{\tilde{B}}_T^{(j)} \left ( \frac{\lfloor t_j N_j \rfloor}{N_j} \right ) \right )^2 \right \vert \\ & \stackrel{\hidewidth \eqref{distancemaxbound} \hidewidth}{\leq} \, \, \, \sup \limits_{\vect \in [0,1]^K} \left \vert  \sum_{j=1}^K  \left ( \frac{1}{\alpha_{T,j}} \tilde{\Delta}_T^{(j)} \left ( t_j  \right )  \right )^2  - \sum_{j=1}^K  \left ( \overline{\tilde{B}}_T^{(j)} \left ( \frac{\lfloor t_j N_j \rfloor}{N_j} \right ) \right )^2 \right \vert  \, \, \stackrel{\hidewidth \eqref{squaredapproxst2} \hidewidth}{=} \, \, \, o_P(1) \end{align*}
$\tilde{P}$-a.s., as $T \to \infty$. This concludes the proof.
$_\Box$\\
\vspace{0.3cm}
\\\textit{Proof of Theorem \ref{squaredapproxtheoremonoriginalprobabilityspace}:}\\
%\vspace{0.1cm}
\\ By Theorem \ref{squaredtheorem} there exist equivalent processes $\{ \tilde{\vecY}_{T,j,i} \, : \, i \geq 1  \} \stackrel{d}{=} \{ \vecY_{T,j,i} \, : \, i \geq 1 \}$ and $K$ Brownian motions $\{\tilde{B}_{T}^{(j)}(t_j) \, : \, t_j \geq 0 \}$ on richer probability spaces $(\tilde{\Omega}^{(j)},\tilde{\mathcal{F}}^{(j)},\tilde{P}^{(j)}), j=1,\dots,K$, such that on the product space $(\tilde{\Omega},\tilde{\mathcal{F}},\tilde{P})$ of these probability spaces
\begin{align} \label{weakapprox1} & \sup \limits_{\vect \in [0,1]^K} \left \vert    \sum_{j=1}^K  \left ( \tilde{\mathcal{D}}_{T}^{(j)}(t_j) \right )^2 - \sum_{j=1}^K \left ( \alpha_{T,j} \tilde{B}_T^{(j)} \left ( \frac{\lfloor t_j N_j \rfloor}{N_j}  \right ) \right )^2  \right \vert  \\ \nonumber & =  \sup \limits_{\vect \in [0,1]^K} \left \vert \sum_{j=1}^K  \left ( \sqrt{N_j} \vecv_{T,j}' (\tilde{\bm \Sigma}_T^{(j)}(t_j) - \bm \Sigma_T^{(j)} (t_j)) \vecw_{T,j} \right )^2 - \sum_{j=1}^K \left ( \alpha_{T,j} \tilde{B}_{T}^{(j)} \left ( \frac{\lfloor t_j N_j \rfloor}{N_j}  \right ) \right )^2 \right \vert = o_P(1),   \end{align}
as $ T \to \infty$, where $\tilde{\bm \Sigma}_T^{(j)}(t_j) = N_j^{-1} \sum_{i=1}^{\lfloor t_j N_j \rfloor} \tilde{\vecY}_{T,j,i} \tilde{\vecY}_{T,j,i}'$ and $\bm \Sigma_{T}^{(j)}(t_j) = E(\tilde{\bm \Sigma}_T^{(j)}(t_j))$ for $t_j \in [0,1], j=1,\dots,K$. Now, since the Skorohod space $( (D[0,1])^K)^\N$ is an uncountable, separable, complete and metric space and thus Borel isomorphic to $[0,1]$ (cf. \cite{Billingsley1999},p.212), i.e. there exists a bijective function $\varphi: ( (D[0,1])^K )^\N \mapsto [0,1]$ such that $\varphi$ and $\varphi^{-1}$ are measurable, and due to the additional assumptions of the Theorem, we can apply \cite[Sec.~21, Lemma~2]{Billingsley1999} with
\begin{align*} \nu = \mathcal{L} \left ( \left \{ \left( \sqrt{N_j} \vecv_{T,j}' (\tilde{\bm \Sigma}_T^{(j)}(\bullet) - \bm \Sigma_T^{(j)} (\bullet)) \vecw_{T,j} \right)_{j=1}^K  : T \ge 1 \right \} , \left \{ \left ( \alpha_{T,j} \tilde{B}_{T}^{(j)} \left ( \frac{\lfloor \bullet N_j \rfloor}{N_j}  \right ) \right )_{j=1}^K : T \ge 1 \right \}  \right )  \end{align*}
and $\sigma = \left \{ \left ( \sqrt{N_j} \vecv_{T,j}' (\hat{\bm \Sigma}_T^{(j)}(\bullet) - \bm \Sigma_T^{(j)} (\bullet)) \vecw_{T,j} \right )_{j=1}^K : T \ge 1 \right \}$ to conclude the existence of \\ $\tau = \left \{  \left ( \alpha_{T,j} B_{T}^{(j)} \left ( \frac{\lfloor \bullet N_j \rfloor}{N_j}  \right )  \right )_{j=1}^K : T \ge 1  \right \} $
as a function of $\sigma $ and $U$, such that
$( \{ B_{T}^{(j)}(t_j) \, : \, t_j \geq 0  \} )_{j=1}^K \stackrel{d}{=} ( \{ \tilde{B}_{T}^{(j)} (t_j) \, : \, t_j \geq 0 \} )_{j=1}^K , $
where $\{ B_{T}^{(j)}(t_j) \, : \, t_j \geq 0 \}$ for $j=1,\dots,K$ are independent Brownian motions defined on the original probability spaces $(\Omega^{(j)},\mathcal{F}^{(j)},P^{(j)})$,
and $ \mathcal{L}( \sigma, \tau ) = \nu $. This implies that
\begin{align*} & \left( \left \{ \sum_{j=1}^K  \left ( \sqrt{N_j} \vecv_{T,j}' (\hat{\bm \Sigma}_T^{(j)}(\bullet) - \bm \Sigma_T^{(j)} (\bullet)) \vecw_{T,j} \right )^2 \right \}_{T \ge 1}  ,  \left \{ \sum_{j=1}^K \left ( \alpha_{T,j} B_{T}^{(j)}  \left ( \frac{\lfloor \bullet N_j \rfloor}{N_j}  \right )  \right )^2  \right \}_{T \ge 1}     \right )
\end{align*}
and 
\begin{align*}
& = \left ( \left \{ \sum_{j=1}^K  \left ( \sqrt{N_j} \vecv_{T,j}' (\tilde{\bm \Sigma}_T^{(j)}(\bullet) - \bm \Sigma_T^{(j)} (\bullet)) \vecw_{T,j} \right )^2  \right \}_{T \ge 1}  , \left \{ \sum_{j=1}^K \left ( \alpha_{T,j} \tilde{B}_{T}^{(j)} \left ( \frac{\lfloor \bullet N_j \rfloor}{N_j}  \right ) \right )^2 \right \}_{T \ge 1}   \right ), \end{align*}
are equal in distribution. Therefore, we may conclude that
\begin{align*}  & \sup \limits_{\vect \in [0,1]^K} \left \vert \sum_{j=1}^K  \left ( \sqrt{N_j} \vecv_{T,j}' (\hat{\bm \Sigma}_T^{(j)}(t_j) - \bm \Sigma_T^{(j)} (t_j)) \vecw_{T,j} \right )^2 - \sum_{j=1}^K \left ( \alpha_{T,j} B_{T}^{(j)} \left ( \frac{\lfloor t_j N_j \rfloor}{N_j}  \right ) \right )^2 \right \vert \\ &  \stackrel{d}{=}  \sup \limits_{\vect \in [0,1]^K} \left \vert \sum_{j=1}^K  \left ( \sqrt{N_j} \vecv_{T,j}' (\tilde{\bm \Sigma}_T^{(j)}(t_j) - \bm \Sigma_T^{(j)} (t_j)) \vecw_{T,j} \right )^2 - \sum_{j=1}^K \left ( \alpha_{T,j} \tilde{B}_{T}^{(j)} \left ( \frac{\lfloor t_j N_j \rfloor}{N_j}  \right ) \right )^2 \right \vert,   \end{align*}
Consequently,  \eqref{weakapprox1} implies 
\begin{align*} 
%&  \sup \limits_{\vect \in [0,1]^K} \left \vert \sum_{j=1}^K  \left ( \mathcal{D}_T^{(j)}(t_j) \right )^2 - \sum_{j=1}^K \left ( \alpha_{T,j} B_{T}^{(j)} \left ( \frac{\lfloor t_j N_j \rfloor}{N_j}  \right ) \right )^2 \right \vert  \\ & = 
\sup \limits_{\vect \in [0,1]^K} \left \vert \sum_{j=1}^K  \left ( \sqrt{N_j} \vecv_{T,j}' (\hat{\bm \Sigma}_T^{(j)}(t_j) - \bm \Sigma_T^{(j)} (t_j)) \vecw_{T,j} \right )^2 - \sum_{j=1}^K \left ( \alpha_{T,j} B_{T}^{(j)} \left ( \frac{\lfloor t_j N_j \rfloor}{N_j}  \right ) \right )^2 \right \vert  = o_P(1)  \end{align*}
as $T \to \infty$. Further, \eqref{weakapproxonoriginal2} follows using the same steps as in the proof of Theorem \ref{squaredtheorem}. $\hfill \Box$\\
\vspace{0.1cm}
\\\textit{Proof of Theorem \ref{Pooledtheorem}:}\\
%\vspace{0.1cm}
\\Let  $\{ \tilde{B}_T^{(j)}(t_j) \, : \, t_j \geq 0 \}$, $j=1,\dots,K$, be the $K$ independent Brownian motions constructed in the proof of Theorem \ref{squaredtheorem} and denote by $\tilde{\mathcal{G}}_T $ the associated version of $ \mathcal{G}_T $. It is easy to verify that  $\{ \tilde{\mathcal{G}}_T(t_1,\dots,t_K) \, : \, t_1,\dots,t_K \geq 0\}$ is a mean zero Gaussian process with finite dimensional marginals having independent increments $ \tilde{\mathcal{G}}_T(t_{1,r},\dots,t_{K,r}) - \tilde{\mathcal{G}}_T (t_{1,r-1},\dots,t_{K,r-1}) $
$t_{j,0},t_{j,1},\dots,t_{j,n} \geq 0, $ $ j=1,\dots,K, n \in \mathbb N,$  $0 = t_{j,0} < t_{j,1} < \dots < t_{j,n}$, $j=1,\dots,K$, given by $  \sum_{j=1}^K \alpha_{T,j} \left [  \tilde{B}_T^{(j)}(t_{j,r}) - \tilde{B}_T^{(j)}(t_{j,r-1})  \right ]  $ and therefore satisfying
\begin{align*} \tilde{\mathcal{G}}_T(t_{1,r},\dots,t_{K,r}) - \tilde{\mathcal{G}}_T (t_{1,r-1},\dots,t_{K,r-1}) & \sim \mathcal{N} \left ( 0 , \sum_{j=1}^K \alpha_{T,j}^2 (t_{j,r} - t_{j,r-1})   \right )    \end{align*} 
It follows that the vector of increments,  
\begin{align*} & \tilde{\bm{\mathrel{G}}}_T(\{ t_{j,r}, j=1,\dots,K,\, r =1,\dots,n   \}) \\  & = \left ( \tilde{\mathcal{G}}_T(t_{1,1},\dots,t_{K,1}) - \underbrace{\tilde{\mathcal{G}}_T(t_{1,0},\dots,u_{K,0})}_{=0} , \dots, \tilde{\mathcal{G}}_T(t_{1,n},\dots,t_{K,n}) - \tilde{\mathcal{G}}_T(t_{1,n-1},\dots, t_{K,n-1}) \right )' \\ & = \left ( \tilde{\mathcal{G}}_T(t_{1,1},\dots,t_{K,1}), \dots, \tilde{\mathcal{G}}_T(t_{1,n},\dots,t_{K,n}) - \tilde{\mathcal{G}}_T(t_{1,n-1},\dots, t_{K,n-1}) \right )',   \end{align*}
is $\mathcal{N}_n(\vecnull,\matA_T)$ distributed, where
\begin{align*} \matA_T = \diag \left (\sum_{j=1}^K \alpha_{T,j}^2 t_{j,1},  \sum_{j=1}^K \alpha_{T,j}^2 (t_{j,2} - t_{j,1} ) , \dots, \sum_{j=1}^K \alpha_{T,j}^2 (t_{j,n} - t_{j,n-1}) \right ) \in \mathbb R^{n \times n}.  \end{align*}
Now let $ \matM = ( \eins(i\le j) )_{1\le i \le n \atop 1 \le j \le n} $. Then 
\begin{align*} & \left ( \tilde{\mathcal{G}}_T (t_{1,1},\dots,t_{K,1}), \tilde{\mathcal{G}}_T (t_{1,2},\dots,t_{K,2}) \dots, \tilde{\mathcal{G}}_T (t_{1,n}, \dots, t_{K,n})   \right )'  \\ & \qquad = \matM \tilde{\bm{\mathrel{G}}}_T (\{ t_{j,r}, j=1,\dots,K, r=1,\dots,n \})' \sim \mathcal{N}_n (\vecnull, \matM \matA_T \matM'), \end{align*}
where
\begin{align*} \matM \matA_\matT \matM' = \begin{pmatrix} \sum_{j=1}^K \alpha_{T,j}^2 t_{j,1} & \sum_{j=1}^K \alpha_{T,j}^2 t_{j,1} & \cdots & \cdots & \cdots &  \sum_{j=1}^K \alpha_{T,j}^2 t_{j,1} \\ \sum_{j=1}^K \alpha_{T,j}^2 t_{j,1} &  \sum_{j=1}^K \alpha_{T,j}^2 t_{j,2} & \sum_{j=1}^K \alpha_{T,j}^2 t_{j,2} & \cdots & \cdots & \sum_{j=1}^K \alpha_{T,j}^2 t_{j,2} \\ \sum_{j=1}^K \alpha_{T,j}^2 t_{j,1} & \sum_{j=1}^K \alpha_{T,j}^2 t_{j,2} & \sum_{j=1}^K \alpha_{T,j}^2 t_{j,3}   & \cdots & \cdots & \sum_{j=1}^K \alpha_{T,j}^2 t_{j,3}      \\ \vdots & \vdots & \vdots & \vdots & \vdots & \vdots \\ \sum_{j=1}^K \alpha_{T,j}^2 t_{j,1} & \sum_{j=1}^K \alpha_{T,j}^2 t_{j,2} & \sum_{j=1}^K \alpha_{T,j}^2 t_{j,3} & \cdots & \cdots & \sum_{j=1}^K \alpha_{T,j}^2 t_{j,n}     \end{pmatrix}. \end{align*}
In particular, the covariance function of $\tilde{\mathcal{G}}_T$ is given
\begin{align*} \gamma_{\tilde{\mathcal{G}}_T}(\vecs,\vect) & = \Cov  ( \tilde{\mathcal{G}}_T (s_1,\dots,s_K)    , \tilde{\mathcal{G}}_T(t_1,\dots,t_K) )  =  \sum_{j=1}^K \alpha_{T,j}^2 \min \{s_j, t_j\}      \end{align*}
for $\vecs = (s_1,\dots,s_K)', \vect = (t_1,\dots,t_K)'$ with $s_j,t_j \geq 0$ for $j=1,\dots,K$.\\
\vspace{0.1cm}
\\Since the partial sum $\tilde{\matS}_{T,k_1,\dots,k_K}$ based on the pooled sample variance-covariance matrix at indices $k_1,\dots,k_K \in \mathbb N_0$ admits the representation $\tilde{\matS}_{T,k_1,\dots,k_K} = \sum_{j=1}^K \hat{\tilde{\bm{\Sigma}}}_{T,k_j}^{(j)}$
\begin{align*} \tilde{D}_{T,k_1,\dots,k_K} = \sum_{j=1}^K \vecv_T' \left ( \hat{\tilde{\bm{\Sigma}}}_{T,k_j}^{(j)} - \bm{\Sigma}_{T,k_j}^{(j)}  \right ) \vecw_T = \sum_{j=1}^K \tilde{D}_{T,k_j}^{(j)}  \end{align*}
and therefore
\begin{align*}  \left \vert \tilde{D}_{T,k_1,\dots,k_K} - \tilde{\mathcal{G}}_{T}(k_1,\dots,k_K)  \right \vert \leq   \sum_{j=1}^K \left \vert \tilde{D}_{T,k_j}^{(j)} - \alpha_{T,j} \tilde{B}_{T}^{(j)}(k_j) \right \vert \stackrel{ \eqref{DTjtapproxx}}{\leq} \sum_{j=1}^K C_{T,j}k_j^{\frac{1}{2} - \lambda_j}, \quad \tilde{P}-\text{a.s.}.  \end{align*}
Similarly we have for the associated \cadlag - version
\begin{align*} & \tilde{\mathcal{D}}_{T}(t_1,\dots,t_K) = \sum_{j=1}^K \frac{1}{\sqrt{N}} \vecv_T' (\hat{\tilde{\bm{\Sigma}}}_{T,\lfloor t_j N_j \rfloor}^{(j)} - \bm{\Sigma}_{T,\lfloor t_j N_j \rfloor}^{(j)} ) \vecw_T = \sum_{j=1}^K \sqrt{\frac{N_j}{N}} \tilde{\mathcal{D}}_{T}^{(j)}(t_j), \quad t_1,\dots,t_K \in [0,1], \end{align*}
and thus receive the strong approximation
\begin{align*} & \sup \limits_{t_1,\dots,t_K \in [0,1]} \left \vert\tilde{\mathcal{D}}_T(t_1,\dots,t_K) -  \tilde{\mathcal{G}}_{T} \left ( \frac{\lfloor  t_1 N_1 \rfloor}{N} , \dots, \frac{\lfloor t_K N_K \rfloor}{N} \right )  \right \vert \\ & \stackrel{\hidewidth \hidewidth}{\leq} \,  \sum_{j=1}^K \sqrt{\frac{N_j}{N}} \sup \limits_{t_j \in [0,1]} \left \vert \tilde{\mathcal{D}}_{T}^{(j)}(t_j) - \alpha_{T,j} \sqrt{\frac{N}{N_j}} \tilde{B}_{T}^{(j)} \left ( \frac{\lfloor t_j N_j  \rfloor}{N} \right ) \right \vert \\ & \stackrel{\hidewidth d \hidewidth}{=} \,   \sum_{j=1}^K \sqrt{\frac{N_j}{N}} \sup \limits_{t_j \in [0,1]} \left \vert \tilde{\mathcal{D}}_{T}^{(j)}(t_j) - \alpha_{T,j} \tilde{B}_{T}^{(j)} \left ( \frac{\lfloor t_j N_j  \rfloor}{N_j} \right ) \right \vert \, \,  \stackrel{\hidewidth \eqref{cadlagapproxx} \hidewidth}{=} \, \,  o(1) \end{align*}
$\tilde{P}$-a.s., for $T \to \infty$ since  $C_{T,j} N_j^{-\lambda_j} = o(1)$ and $\frac{N_j}{N} \to \kappa_j \in (0,1]$ for all $j=1,\dots,K$ as $T \to \infty$, by assumption. $ \hfill \Box$\\
\vspace{0.2cm}
\\\textit{Proof of Corollary \ref{pooledcorollary1}:}\\
%\vspace{0.1cm}
\\(i) Using \eqref{pooledapproxx} we have 
\begin{align*} & \sup \limits_{ \vect \in [0,1]^K}  \left \vert  \tilde{\mathcal{D}}_T^* (t_1,\dots,t_K) - \tilde{\mathcal{G}}_T^* \left ( \frac{\lfloor t_1 N_1 \rfloor}{N}, \dots, \frac{\lfloor t_K N_K \rfloor}{N} \right )   \right \vert  \\ & = \sup \limits_{\vect \in [0,1]^K} \left \vert \frac{1}{\sqrt{\sum_{j=1}^K \alpha_{T,j}^2}} \tilde{\mathcal{D}}_T(t_1,\dots,t_K) -  \frac{1}{\sqrt{\sum_{j'=1}^K \alpha_{T,j'}^2}} \tilde{\mathcal{G}}_T \left ( \frac{\lfloor t_1 N_1 \rfloor}{N}, \dots, \frac{\lfloor t_K N_K \rfloor}{N} \right )  \right \vert \\ & = \frac{1}{\sqrt{\sum_{j'=1}^K \alpha_{T,j'}^2}} \sup \limits_{\vect \in [0,1]^K} \left  \vert \tilde{\mathcal{D}}_T(t_1,\dots,t_K) - \tilde{\mathcal{G}}_T \left ( \frac{\lfloor t_1 N_1 \rfloor}{N}, \dots, \frac{\lfloor t_K N_K \rfloor}{N} \right ) \right \vert \\ &  \stackrel{\hidewidth \eqref{pooledapproxx} \hidewidth}{=} \, \, \, o(1)    \end{align*}
$\tilde{P}$-a.s., as $T \to \infty$ since $\alpha_{T,j} \to \alpha_j^* \in (0,\infty)$ as $T \to \infty$.\\
\vspace{0.2cm}
\\(ii) The result follows, if we show the CLT for the equivalent statistic $ \tilde{S}_T^* $, since then $ P( S_T^* \le x ) = \wt{P}( \tilde{S}_T^* \le x ) \to \Phi(x) $, as $ T \to \infty $, for any $ x \in \mathbb{R} $. Denote the equivalent version of $\{\mathcal{D}_T(t_1,\dots,t_K) \, : \, t_1,\dots,t_K  \in [0,1]\}$, which appears in Theorem \ref{Pooledtheorem} and is defined on the probability space $(\tilde{\Omega},\tilde{F},\tilde{P}),$ by $\{ \tilde{\mathcal{D}}_T(t_1,\dots,t_K) \, : \, t_1,\dots,t_K \in [0,1]\}$. Further let us define 
\begin{align*} \tilde{\mathcal{D}}_T^*(t_1,\dots,t_K)  = \sqrt{  \frac{N}{\sum_{j'=1}^K \alpha_{T,j'}^2 N_{j'}}} \tilde{\mathcal{D}}_T(t_1,\dots,t_K), \quad t_1,\dots,t_K \in [0,1], \end{align*}
and the corresponding scaled version of $ \tilde{\mathcal{G}}_T $, $\tilde{\mathcal{G}}_T^*(u_1,\dots,u_K) = \sqrt{N \slash \sum_{j'=1}^K \alpha_{T,j'}^2 N_{j'}}  \tilde{\mathcal{G}}_T \left ( u_1,\dots,u_K \right )$ for \\ $u_1,\dots,u_K \geq 0$. The assertion follows at once by noting that
\begin{align*} \tilde{\mathcal{D}}_T^*(1,\dots,1) = \frac{N}{\sqrt{\sum_{j'=1}^K \alpha_{T,j'}^2 N_{j'}}} \vecv_T' (\tilde{\matS}_{T} - E(\tilde{\matS}_{T})) \vecw_T ,  \end{align*}
and observing that 
\begin{align*} \tilde{\mathcal{G}}_T^* \left ( \frac{N_1}{N}, \dots, \frac{N_K}{N}  \right ) &  = \sqrt{\frac{N}{\sum_{j'=1}^K \alpha_{T,j'}^2 N_{j'}}} \sum_{j=1}^K \alpha_{T,j} \left [ \tilde{B}_T^{(j)} \left ( \frac{N_j}{N}  \right )  - \underbrace{\tilde{B}_T^{(j)}(0)}_{=0} \right ] \\ & \sim \sqrt{\frac{N}{\sum_{j'=1}^K \alpha_{T,j'}^2 N_{j'}}} \sum_{j=1}^K \alpha_{T,j} \tilde{X}_j, \qquad \tilde{X}_j \sim \mathcal{N} \left (0, \frac{N_j}{N} \right ), \\ & \sim \mathcal{N}(0,1)  \end{align*}
for all $T > 0$. $ \hfill \Box$\\
\vspace{0.2cm}
\\\textit{Proof of Theorem \ref{Thmdtschlangeapprox}:}\\ 
%\vspace{0.1cm} 
\\We have $\sup_{t_1,\dots,t_K \in [0,1]} \left \vert \tilde{\Delta}_{T}(t_1,\dots,t_K) - \tilde{\Gamma}_{T} \left ( t_1,\dots,t_K   \right ) \right \vert \leq I_1 + I_2 $
with 
\begin{align*} I_1 = \sup \limits_{t_1,\dots,t_K \in [0,1]} \left \vert \tilde{\mathcal{D}}_{T} \left ( t_1,\dots,t_K  \right ) -  \tilde{\mathcal{G}}_{T} \left ( \frac{\lfloor t_1 N_1 \rfloor}{N}, \dots, \frac{\lfloor t_K N_K \rfloor}{N}   \right ) \right \vert  \stackrel{\eqref{pooledapproxx}}{=} o(1),  \end{align*}
\begin{align*} I_2 &= \sup \limits_{t_1,\dots,t_K \in [0,1]} \left \vert  \sum_{j=1}^K \frac{\lfloor t_j N_j \rfloor}{\sqrt{N}} \frac{1}{\sqrt{N_j}} \tilde{\mathcal{D}}_{T}^{(j)}(1)  -  \sum_{j=1}^K \frac{\lfloor t_j N_j \rfloor}{\sqrt{N}} \frac{1}{\sqrt{N_j}} \alpha_{T,j} \tilde{B}_{T}^{(j)}(1)  \right \vert \\ & \leq \sum_{j=1}^K   \sqrt{\frac{N_j}{N}} \sup \limits_{t_j \in [0,1]} \underbrace{\frac{\lfloor t_j N_j \rfloor}{N_j}}_{\leq 1}  \left \vert \tilde{\mathcal{D}}_{T}^{(j)}(1) - \alpha_{T,j} \tilde{B}_{T}^{(j)}(1)     \right \vert  \stackrel{\eqref{cadlagapproxx}}{=} o(1),  \end{align*}
as $T \to \infty$, when $\frac{N_j}{N} \to \kappa_j \in (0,1]$ for $j=1,\dots,K$ as $T \to \infty$. Further, the Brownian bridge representation of $\tilde{\Gamma}_T$ follows from the scaling property of Brownian motions,
\begin{align*} \tilde{\Gamma}_T(t_1,\dots,t_K) &= \sum_{j=1}^K \alpha_{T,j} \tilde{B}_T^{(j)} \left ( \frac{\lfloor t_j N_j \rfloor}{N}  \right ) - \frac{1}{\sqrt{N}} \sum_{j=1}^K \frac{\lfloor t_j N_j \rfloor}{\sqrt{N_j}} \alpha_{T,j} \tilde{B}_T^{(j)}(1) \\ & \stackrel{d}{=} \sum_{j=1}^K \alpha_{T,j} \sqrt{\frac{N_j}{N}} \left [ \tilde{B}_T^{(j)} \left ( \frac{\lfloor t_j N_j \rfloor}{N_j}  \right ) - \frac{\lfloor t_j N_j \rfloor}{N_j} \tilde{B}_T^{(j)}(1)   \right ], \end{align*}
which establishes the result. $ \hfill \Box$\\
\vspace{0.2cm}
\\\textit{Proof of Corollary \ref{cusumcorollary2}: }
The proof is straightforward and omitted. $ \hfill \Box $
\vspace{0.3cm}
\\\textit{Proof of Theorem \ref{pooledapproxonoriginalspace}:} The proof can be found in \cite{Mause2019}. $\hfill \Box$
\vspace{0.3cm}
\\\textit{Proof of Lemma \ref{teststatisticslemma}:}\\
\vspace{0.1cm}
\\We show (ii). The other assertions follow along the same lines. First note that there exists, for all $j = 1,\dots,K $, on $(\Omega^{(j)},\mathcal{F}^{(j)},P^{(j)})$, a standard Brownian motion $\{ B^{(j)}(t_j) \, : \, t_j \geq 0\}$, such that for all $T > 0$ $ \{ B^{(j)}(t_j) \, : \, t_j \geq 0 \} \stackrel{d}{=} \{ B_T^{(j)}(t_j) \, : \, t_j \geq 0 \}, $ % \end{align*}
and thus, by independence of the $K$ Brownian motions, $\bm B_T(\vect) \stackrel{d}{=} \bm B(\vect)$. Similarly one can show $ \overline{ \bm B}_T(\vect) \stackrel{d}{=}  \overline{\bm B}(\vect)$ for the vector of Brownian bridges $\overline{\bm B}_T(\vect) = ( \overline{B}_T^{(1)}(t_1), \dots, \overline{B}_T^{(K)}(t_K ))'$.\\
\vspace{0.1cm}
We have
\begin{align*} & \left \vert V_T -  \max \limits_{k \in \mathcal{M}_T} \left \vert \mathcal{G}_T \left ( \frac{k_1}{N},\dots,\frac{k_K}{N}  \right ) \right \vert        \right \vert  =  \left \vert V_T -  \max \limits_{k \in \mathcal{M}_T} \left \vert \sum_{j=1}^K \alpha_{T,j} B_T^{(j)} \left ( \frac{k_j}{N}  \right ) \right \vert        \right \vert     = o_P(1)  \end{align*}
as $T \to \infty$ on $(\Omega,\mathcal{F},P)$, where (using the scaling property of Brownian motion)
\begin{align*} \max \limits_{\veck \in \mathcal{M}_T} \left \vert  \sum_{j=1}^K \alpha_{T,j} B_T^{(j)} \left (\frac{k_j}{N}  \right )  \right \vert \stackrel{d}{=} \max \limits_{\veck \in \mathcal{M}_T} \left \vert  \sum_{j=1}^K \alpha_{T,j} \sqrt{\frac{N_j}{N}} B^{(j)} \left (\frac{k_j}{N_j}  \right ) \right \vert  
\end{align*}
Therefore, 
\begin{align*} \left \vert  V_T - \max \limits_{\veck \in \mathcal{M}_T} \left \vert  \sum_{j=1}^K \alpha_{T,j}  B_T^{(j)} \left (\frac{k_j}{N}  \right )  \right \vert  \right \vert
\stackrel{d}{=} \left \vert V_T -  \sup \limits_{\vect \in  [0,1]^K } \left \vert  \sum_{j=1}^K \alpha_{T,j} \sqrt{\frac{N_j}{N}} B^{(j)} \left ( \frac{\lfloor t_j N_j \rfloor}{N} \right )   \right \vert   \right  \vert \end{align*}
and consequently, $\left \vert   V_T - \sup \limits_{\vect \in  [0,1]^K } \left \vert  \sum_{j=1}^K \alpha_{T,j} \sqrt{\frac{N_j}{N}} B^{(j)} \left ( \frac{\lfloor t_j N_j \rfloor}{N} \right )   \right \vert   \right \vert = o_P(1),\, T \to \infty. $ Since the paths of $B^{(j)}$ are a.s. continuous it follows that 
\begin{align*} \sup \limits_{\vect \in  [0,1]^K } \left \vert  \sum_{j=1}^K \alpha_{T,j} \sqrt{\frac{N_j}{N}} B^{(j)} \left ( \frac{\lfloor t_j N_j \rfloor}{N_j} \right )   \right \vert \longrightarrow \sup \limits_{\vecs \in  [0,1]^K } \left \vert  \sum_{j=1}^K \alpha_{j}^* \sqrt{\kappa_j} B^{(j)} \left ( s_j \right )   \right \vert , \quad T \to \infty,   \end{align*}
$P$-a.s., and hence in distribution, where the limiting r.v. is a.s. finite, since $\alpha_{T,j} \to \alpha_j^*$ and $\frac{N_j}{N} \to \kappa_j$ as $T \to \infty$ for $j=1,\dots,K$.
Consequently, we can conclude that 
%\begin{align*} 
$V_T \stackrel{d}{\longrightarrow}  \sup \limits_{\vecs \in  [0,1]^K } \left \vert  \sum_{j=1}^K \alpha_{j}^* \sqrt{\kappa_j} B^{(j)} \left ( s_j \right )   \right \vert , \quad T \to \infty, $ %\end{align*}
by Theorem 2.7, (iv) in \cite{vandervaart1998}.$\hfill _\Box$

\section*{ACKNOWLEDGEMENTS}
The authors gratefully acknowledge the financial support from Deutsche Forschungsgemeinschaft (DFG), grants  STE 1034/11-1 and 1034/11-2, and would like to thank the anonymous referees for their careful reading of the manuscript.

\end{document}